\documentclass[12pt]{article}

\usepackage{mathtools,amsfonts,amssymb}
\usepackage{array}
\usepackage{epsfig}
\usepackage{graphicx}
\usepackage[centertableaux,boxsize=1.1em]{ytableau}
\usepackage{hyperref}
\usepackage{bm}
\usepackage{fancyhdr}
\usepackage[margin=1in]{geometry}
\usepackage{booktabs}
\usepackage[shortlabels,inline]{enumitem}
\usepackage[nameinlink,capitalize]{cleveref}

\newtheorem{conj}{Conjecture}
\newtheorem{cor}{Corollary}
\newtheorem{prop}{Proposition}
\newtheorem{thm}{Theorem}

\DeclareMathOperator{\SU}{SU}
\DeclareMathOperator{\U}{U}

\renewcommand{\bar}[1]{\ensuremath\overline{#1}}

\newcommand\bC{\mathbb{C}}
\newcommand\bD{\mathbb{D}}

\newcommand\bH{\mathbb{H}}
\newcommand\bL{\mathbb{L}}
\newcommand\bP{\mathbb{P}}
\newcommand\bQ{\mathbb{Q}}
\newcommand\bR{\mathbb{R}}
\newcommand\bV{\mathbb{V}}
\newcommand\bZ{\mathbb{Z}}

\newcommand\B{\mathcal{B}}

\newcommand\cD{\mathcal{D}}
\newcommand\cF{\mathcal{F}}
\newcommand\cH{\mathcal{H}}

\newcommand\cJ{\mathcal{J}}
\newcommand\cM{\mathcal{M}}
\newcommand\cO{\mathcal{O}}
\newcommand\cP{\mathcal{P}}

\newcommand\SL{\operatorname{SL}}

\makeatletter
\def\blfootnote{\gdef\@thefnmark{}\@footnotetext}
\makeatother

\title{Dimensional Reduction of B-fields in F-theory}
\author{Sheldon Katz\\[0.9em]
\small\textit{Department of Mathematics} \\
\small\textit{University of Illinois at Urbana-Champaign} \\%
\small\textit{1409W.\ Green St.} \\%
\small\textit{Urbana, IL 61801, USA}\\[1em]
Washington Taylor\\[0.9em]
\small\textit{Center for Theoretical Physics, Department of Physics} \\
\small\textit{Massachusetts Institute of Technology} \\%
\small\textit{77 Massachusetts Avenue} \\%
\small\textit{Cambridge, MA 02139, USA}%
}
\date{\today}

\begin{document}
\maketitle
\thispagestyle{fancy}

\vspace{-2em}
\begin{abstract}
We describe the dimensional reduction of the IIB B-fields in F-theory using a conjectured description of normalizable B-fields in terms of perverse sheaves.  Computations are facilitated using the Decomposition Theorem.  Many of our descriptions are new, and all our results are all consistent with known results in physics.  We also conjecture a physical framework for normalizable B-fields and show consistency with mathematics. 

\smallskip
We dedicate this paper to Herb Clemens, in admiration for his myriad fundamental contributions to complex algebraic geometry, together with his more recent  interest in F-theory in physics.  This paper deals with three of Herb's interests: Hodge theory, topology of algebraic varieties, and F-theory, and so is a fitting way for us to express our appreciation for his contributions over a period of more than five decades.
\end{abstract}

\blfootnote{
\texttt{wati} at \texttt{mit.edu},
\texttt{katzs} at \texttt{illinois.edu}
}
\vspace{-1.2em}

\newpage
\tableofcontents

\newpage

\section{Introduction}\label{sec:intro}



In this paper, we consider F-theory as type IIB string theory with varying axio-dilaton, associated with elliptically fibered Calabi-Yau's $\pi:X\to \B$.  For the benefit of mathematicians not familiar with string theory or F-theory, we unpack some of this terminology below, which should suffice for understanding the mathematical content of our work and how it relates to physics.  The papers \cite{Vafa-F-theory,Morrison-Vafa-I,Morrison-Vafa-II} explain some of the basic ideas of F-theory and how it relates to algebraic geometry.  A version of our work oriented towards physicists will appear elsewhere \cite{kt}; a brief summary of this work from the
physics perspective is given in \S\ref{sec:physics}.

Our starting point is the work of \cite{dps}, which considered the
case of a generic elliptically fibered K3 surface $\pi:S\to \bP^1$,
with nodal fibers over 24 points of $\bP^1$.  That paper resolved a
long-standing puzzle based on the discrepancy between the $U(1)^{20}$
gauge group expected for a generic elliptically fibered K3 by duality
with M-theory, and an apparent $U(1)^{24}$ gauge group arising from
7-branes located at the 24 points of the discriminant.  The
discrepancy was resolved using two ideas.  First, the $U(1)^{24}$ is
argued to be nonphysical via the Cremmer-Scherk mechanism \cite{cs}.
Second, it was explained how the missing 20 gauge fields arise from
dimensional reduction of the doublet of B-fields of IIB string theory
after requiring a normalizability condition.  The mathematical
justification of the calculation was based on $L^2$ and
Hodge-theoretic results of Zucker \cite{Zucker}.

In this paper, we extend the analysis of \cite{dps} to elliptically
fibered Calabi-Yau threefolds.  Up until now, the precise field
content of these 6D F-theory compactifications could only be
determined by duality with M-theory.  We show that if we view F-theory
as IIB with varying axio-dilaton without recourse to duality, then
both gauge fields and scalar fields can arise from dimensional
reductions of the B-fields and we conjecture a precise mathematical
framework
for describing these fields.  Our results clarify how F-theory can be
defined directly in terms of type IIB string theory without recourse
to an M-theory limit, and are consistent with
other considerations of physics, including anomaly cancellation in 6
dimensions.
Some other recent works that focus on the related question of defining
F-theory from geometry without resolution are \cite{Jefferson:2021bid,Grassi:2021ptc}.

A full description of the fields of F-theory requires a more detailed
analysis of the 7-branes than we give here, together with an extension
of the application of the Cremmer-Scherk mechanism employed in
\cite{dps}.  
In particular, there are subtleties related to nonabelian gauge fields and
localized scalar fields
associated with charged hypermultiplets in 6D supergravity whose
further elucidation we leave to future work.

\smallskip
We now give a quick overview of IIB string theory and F-theory; these
ideas are expanded on further in \S\ref{sec:physics}.  Type IIB string
theory is a physical theory in 10 (real) spacetime dimensions.  This
theory contains \emph{fields} that permeate the 10-dimensional
spacetime: two scalar fields (the axion and the dilaton), two 2-form
gauge fields called B-fields, and the graviton.  There is a precise
mathematical description of the $p$-form gauge fields appearing in
this paper in terms of Deligne cohomology.  There are also extended
objects of odd spatial dimensions called \emph{branes} that move in
time and play a fundamental role in the theory.  In physical theories,
massless scalar fields parametrize the moduli space of the theory.

In addition, type IIB string theory has an $\SL(2,\bZ)$-symmetry, which can be described geometrically using a principal $SL(2,\bZ)$-bundle $\cP$ on the spacetime.    The axion $\chi$ and dilaton $\phi$ combine to form the complex-valued axio-dilaton $\tau=\chi+ie^{-\phi}$, which transforms in the usual way under $\SL(2,\bZ)$. The pair of B-fields transforms under $\SL(2,\bZ)$ as the usual 2-dimensional representation.   More invariantly, these objects can be described in terms of sections of rank 2 vector bundles on spacetime associated to $\cP$ and the fundamental two-dimensional representation of $\SL(2,\bZ)$.

Another key idea is \emph{compactification}, where a compact space (such as a smooth projective variety) is used to produce an \emph{effective physical theory} of lower dimension by \emph{dimensional reduction} of the constituent fields.   We sketch how this works with F-theory.  

We start with an elliptically fibered Calabi-Yau $n$-fold $\pi:X\to \B$, described by a Weierstrass fibration over a smooth projective $\B$: 
\begin{equation}\label{eq:weierstrass}
X\subset \bP(\cO_B(-2K_B)\oplus  \cO_B(-3K_B)\oplus \cO_B),\qquad y^2z=x^3+fxz^2 +gz^3,
\end{equation}
where $f\in H^0(\B,\cO_B(-4K_B))$, $g\in H^0(\B,\cO_B(-6K_B))$, and $(x,y,z)$ lives in the projective bundle.
We do not need to assume that $X$ is smooth, as the dualizing sheaf of $X$ is trivial by adjunction, irrespective of the singularities of $X$.  The singular fibers of $\pi$ are parametrized by the discriminant divisor
\begin{equation}
\Delta\subset \B, \qquad 4f^3+27g^2=0.
\end{equation}
F-theory also has a moduli space of K\"ahler metrics, as we will elaborate on later.  We put $U=\B-\Delta$.

We now take the 10-dimensional spacetime to be $\B\times M^{12-2n}$, where
$M^{12-2n}$ is $(12-2n)$-dimensional Minkowski space.  The compactified theory will be a physical theory on $M^{12-2n}$.  For an elliptically fibered K3 surface we have $n=2$ and get an 8-dimensional theory.  For an elliptically fibered Calabi-Yau threefold we have $n=3$ and get a 6-dimensional theory.

A configuration of
7-branes\footnote{The 7-brane is being viewed as the
  7-dimensional Riemannian manifold $\Delta\times\bR^{11-2n}$ moving along the
  time direction of $M^{12-2n}$.} is supported on
$\Delta\times M^{12-2n}$.\footnote{More precisely, there can be more
  than one 7-brane, and the 7-branes can be supported on components of
  $\Delta$ rather than $\Delta$ itself.}  The fields of
the 
bulk
supergravity
 theory are generally
not defined along the 7-brane, much as the familiar electric field
blows up at the location of a charged particle.  Thus the supergravity
fields are
defined on $U\times M^{12-2n}$.  
In string theory additional fields are localized on the 7-branes,
described by open strings.  The goal of this work is to elucidate how
the fields living on the 7-branes can be understood in terms of
supergravity fields localized near these branes.
This picture in some sense is an inversion of the holographic picture
in which the fields on the branes capture the near-horizon physics of
the bulk supergravity theory.
In our analysis, we let $\pi_1:U\times M^{12-2n}\to U$ be the
projection.

We now relate the axio-dilaton $\tau$ and the B-fields to the geometry
of the elliptic fibration.  Given $b\in U$, the elliptic curve
$E_b=\pi^{-1}(b)$ can be written as $E_b\simeq \bC/(\bZ+\tau(b)\bZ)$
for some $\tau(b)$ in the upper half plane, well-defined up to
$\SL(2,\bZ)$, so that $\tau$ becomes a multivalued function on $U$,
transforming under $\SL(2,\bZ)$ in the usual way.  Pulling back by
$\pi_1$, we get a multivalued function on $U\times M^{12-2n}$,
also denoted by $\tau$, transforming under $\SL(2,\bZ)$ in
the usual way.  In this way, the part of the moduli of F-theory
parametrized by the axio-dilaton exactly matches the moduli space of
elliptic fibrations $X\to \B$.

Now let $\pi_U:\pi^{-1}(U)\to U$ be the restriction of $\pi$ over $U\subset \B$.  Let $\bV=R^1(\pi_U)_*(\bR)$, the flat rank-2 real vector bundle on $U$ associated to the first real cohomology of the elliptic fibers.  Then the B-fields are closed 2-forms on $U\times M^{12-2n}$ with values in $\bV$.

We next explain dimensional reduction.  Let $\omega$ be a $2-p$-form gauge field on $M^{12-2n}$, and let $\eta$ be a $p$-form on $U$ with values in $\bV$, $0\le p\le 2$.  Then
\begin{equation}\label{eq:reduceb}
B = \pi_1^*(\eta)\wedge \pi_2^*(\omega)
\end{equation}
is a two-form gauge field on $U\times M^{12-2n}$ with values in $\bV$;
i.e.\ B is a B-field.  We say that the $2-p$-form gauge field $\omega$
arises by dimensionally reducing $B$ along $\eta$.  The form $\omega$
has an interpretation in a $(12-2n)$-dimensional physical theory on $M^{12-2n}$.
These fields are parametrized by the $\eta$.  While clearly there are
many more B-fields than those expressible in the form of
(\ref{eq:reduceb}), the B-fields given by (\ref{eq:reduceb}) are the
ones that are relevant for dimensional reduction.
When $p=2$, $\omega$ is a function, interpreted as a
scalar field, a moduli parameter for the $(12-2n)$-dimensional theory.  When
$p=1$, $\omega$ is an ordinary (1-form) gauge field, corresponding to
a gauge symmetry of the theory.

Because we are interested in massless $(2-p)$-form gauge
fields in the dimensionally reduced theory (there are also massive
fields at high energy scales), we are interested in $p$-forms $\eta$
that are closed.  Generalizing the
picture developed in \cite{dps}, we make a physics conjecture that the
forms that are $L^2$ with respect to the natural physical K\"ahler
metric on $U$ correspond to the fields associated with bulk
supergravity degrees of freedom, while those forms that are not
normalizable correspond to degrees of freedom living on the 7-branes.
The corresponding $L^2$ cohomology classes for the bulk fields have
harmonic representatives that minimize the norm.

Our main mathematical conjecture is that the physically relevant $L^2$
cohomologies or spaces of harmonic forms can be identified with
$\bH^p(\B,IC(\bV))$, the $p$th hypercohomology of the intersection
cohomology complex associated to $\bV$.  This conjecture is naturally
suggested by considering the Decomposition Theorem of \cite{bbd} in
comparing F-theory and M-theory, and is consistent with analogous
results in arbitrary dimension in a more general setting where $\bV$
is a polarized variation of Hodge structure on the complement of a
normal crossings divisors \cite{cks2,KK}.  The papers
\cite{Zucker,cks2,KK} were inspired by an unpublished conjecture of
Deligne.  Our conventions are such that $IC(\bV)[\dim \B]$ is a
perverse sheaf on $\B$.

\smallskip
Here is an outline of the paper.  

In Section~\ref{sec:hodge} we review the results of \cite{Zucker,cks2,KK}, which equates the $L^2$ cohomology of a polarized variation of Hodge structure $\bV$ on the complement of a normal crossings divisor with the hypercohomology of $IC(\bV)$ under certain hypothesis on the asymptotic form of the K\"ahler metric.  We also provide some general background on $IC(\bV)$, and perverse sheaves more generally.

In Section \ref{sec:physics} we summarize the physics of 8D and 6D
F-theory compactifications, state the physics conjecture about the
decomposition of forms into normalizable and non-normalizable
associated with bulk and seven-brane degrees of freedom, and give some
examples of specific theories where our conjectures can be tested.

In Section \ref{sec:main-conjecture} we state our main math
conjectures and apply the Decomposition Theory to compute the scalars
and gauge fields obtained by dimensionally reducing the B-fields.  We
check consistency of our results with known properties of F-theory and
find complete agreement.

In Section \ref{sec:MW}, we construct a map from the Mordell-Weil group of an elliptically fibered Calabi-Yau threefold to $H^1(\B,IC(\bV))$.  Given a rational section of the elliptic fibration, we also propose a de Rham representative of the restriction to $U$ of its class in $H^1(\B,IC(\bV))$.  We also check that this de Rham cohomology class vanishes for torsion sections, as expected from physics.

In Section \ref{subsec:supports}, we use cohomology with supports to
illustrate the holographically inspired decomposition of
fields in 8-dimensional F-theory, giving further perspective on the
work of \cite{dps}.  We find a 44-dimensional space of
cohomology classes $H^1(U,\bV)$ with a 20-dimensional subspace
$H^1(\bP^1,j_*\bV)$ which we argue is related to normalizable 1-forms,
and a 24-dimensional quotient space associated with data on the
7-branes.  We extend our proof to 6-dimensional F-theory in the case where all elliptic fibers are of Kodaira type $I_1$, and give a plausability argument for why this might work in complete generality.

\smallskip\noindent \emph{Acknowledgements.} We thank Eduardo Cattani
for helpful conversations and for providing references, especially
\cite{cks2,KK}.  We also thank David Cox, Michael Douglas, Daniel Park
and Lutian Zhao for helpful conversations.  The work of SK is
partially supported by NSF grant DMS-1802242.
The work of WT is supported by
DOE
grant DE-SC00012567.

\section{$L^2$ cohomology and intersection cohomology}\label{sec:hodge}

\subsection{$L^2$ cohomology and intersection cohomology in the Poincar\'e metric}\label{subsec:poincare}

In this section, we review the results of \cite{Zucker,cks2,KK}.  We
let $\B$ be a compact K\"ahler manifold of arbitrary dimension $n$ and
let $D\subset \B$ 
be
a normal crossings divisor.  Let $U=\B-D$ and let
$\bV_{\bC}$\footnote{We use the notation $\bV_{\bC}$ to denote a local
  system of complex vector spaces.  The notation $\bV$ will be
  reserved for local systems of real vector spaces.} be a polarized
variation of Hodge structure on $U$ \cite{griffiths}.  The typical
situation is where $U$ parametrizes a family of smooth projective
varieties $\{X_u\}_{u\in U}$ and the fibers of $\bV_{\bC}$ are
$(\bV_{\bC})_u=H^m(X_u,\bC)$ for some fixed $m$.  As vector spaces,
the fibers $(\bV_{\bC})_u$ are isomorphic to the same vector space $V$
(but can undergo monodromy).  These vector spaces all carry pure Hodge
structures of weight $m$.

Consider a neighborhood $W\subset \B$ of a point $p\in D$ with
$W\simeq \Delta^n$ and $W\cap U\simeq(\Delta^*)^{k}\times
\Delta^{n-k}$ for some $k$.\footnote{$\Delta^n$ is the $n$-disk and $\Delta^*$ is the punctured disk,
while only the unadorned $\Delta$ refers to the discriminant as above.}
We choose a K\"ahler metric on $U$ whose
restriction to $W\cap U$ is asymptotic to
\begin{equation}\label{eq:asympmetric}
    \frac{i}2\left(\sum_{i=1}^{k}\frac{dz_i\wedge d\bar{z}_i}{|z_i|^2\log^2|z_i|^2}+\sum_{i=k+1}^{n} dz_i\wedge d\bar{z}_i\right),
\end{equation}
a product of Euclidean and Poincar\'e metrics.  Such metrics exist \cite{cg}.

We consider the sheaves $L^p_{(2)}(\bV_{\bC})$ on $\B$ of
$\bV_{\bC}$-valued forms which are $L^2$ on compact subsets of $\B$
with locally $L^2$ derivatives, so the sheaves $L^p_{(2)}(\bV_{\bC})$
form a complex $L^\bullet_{(2)}(\bV_{\bC})$.  The norm is defined
using the metric on $\B$ and the Hodge metric on
$\bV_{\bC}$.  It can be shown that the Hodge inner product coincides with
the inner product given by the supergravity theory
 (\ref{eq:m-physics}) \cite{kt}. 

It can be shown that the complex $L^\bullet_{(2)}(\bV_{\bC})$ is well-defined, independent of the choice of K\"ahler metric asymptotic to (\ref{eq:asympmetric}) along $D$.  The sheaves $L^p_{(2)}(\bV_{\bC})$ are fine, so the hypercohomologies $\bH^*(X,L^\bullet_{(2)}(\bV_{\bC}))$ can be computed as the cohomology of the complex of global $L^2$ forms $\Gamma(X,L^\bullet_{(2)}(\bV_{\bC}))$.

In \cite{Zucker}, it was shown 
 for $n=1$, further assuming that
the monodromies of $\bV_\bC$ are unipotent, that this $L^2$ cohomology
is isomorphic to the space of harmonic $\bV_{\bC}$-valued forms on $U$
and also to the sheaf cohomology $H^*(X,j_*\bV_{\bC})$, where $j:U\to
\B$ is the inclusion.  Furthermore, it was shown that
$H^i(X,j_*\bV_{\bC})$ carries a natural Hodge structure of weight
$i+m$, where $m$ is the weight of the variation of Hodge structure
$\bV_{\bC}$.
The normal crossings hypothesis was not needed in the $n=1$ case since in dimension 1, any nonempty Zariski open subset $U\subset \B$ is the complement of finitely points, which is trivially a normal crossings divisor.  

The case of general dimension $n$ with normal crossings boundary was handled in \cite{cks2,KK}.  In this more general situation, the $L^2$ cohomology is again shown to be isomorphic to the space of harmonic $\bV_{\bC}$-valued forms on $U$.  The hypothesis that the monodromy is unipotent was not needed there.

In the case $n>1$, examples can be given where the space of harmonic $p$-forms valued in $\bV_{\bC}$ is \emph{not} isomorphic to $H^p(X,j_*\bV_{\bC})$.  Instead, the space of these harmonic forms is isomorphic to a hypercohomology group to be described presently in terms of perverse sheaves \cite{bbd}.  

\subsection{Perverse Sheaves}\label{subsec:perverse}
Consider the derived category of complexes of sheaves on $\B$ of vector spaces over a fixed field $k$, whose cohomology sheaves $\cH^i(\cF^\bullet)$ are constructible (locally constant on a good stratification of $\B$ and finite dimensional).  We call this category the \emph{constructible derived category}.  Then the 
perverse sheaves $\cF^\bullet$ on $\B$ are objects of the constructible derived category which satisfy a condition on their supports ($\dim \operatorname{supp}(\cH^i(\cF^\bullet))\le -i$ for all $i$) and cosupports ($\dim\operatorname{supp}(\cH^i(\bD\cF^\bullet))\le -i$ for all $i$), where $\bD$ is the Verdier duality functor $\mathrm{RHom}(-,\B[2n])$, and $\operatorname{supp}$ denotes the support of a sheaf.  In this paper we are primarily concerned with $k=\bR$ and $k=\bC$.  Some of the references we provide below use $k=\bQ$ or $k=\bQ_\ell$.  We leave it to the reader to make the necessary adjustments.

As an aside for $k=\bC$, we note that perverse sheaves can alternatively be described as the objects associated to regular holonomic $\cD_B$-modules by the Riemann-Hilbert correspondence, where $\cD_B$ is the sheaf of differential operators on $\B$.  More precisely, the Riemann-Hilbert correspondence is an equivalence of categories between the category of regular holonomic $\cD_B$-modules and the category of perverse sheaves on $\B$.
If $\bV_{\bC}$ is a local system on all of $\B$ then the flat connection endows the vector bundle $\bV_{\bC}\otimes \cO_B$ with a $\cD_B$-module structure.  The Riemann-Hilbert correspondence associates the perverse sheaf $^\pi\bV_{\bC}=\bV_{\bC}[n]$ to this $\cD_B$-module.  This perverse sheaf is represented by the complex whose only nonzero term consists of $\bV_{\bC}$ in degree $-n$. 

For later use, we emphasize that the cohomology sheaves $\cH^i(\cF^\bullet)$ are the usual cohomology sheaves associated to the $\cF^\bullet$, i.e.\ 
\begin{equation*}
  \cH^i(\cF^\bullet)=  \ker\left(d^i:\cF^i\rightarrow \cF^{i+1}\right)/\mathrm{im}\left(d^{i-1}:\cF^{i-1}\rightarrow \cF^{i}\right),
\end{equation*}
and as such they are ordinary sheaves.  Said differently, the $\cH^i(\cF^\bullet)$ are the cohomology sheaves with respect to the standard t-structure, whose heart is the abelian category of constructible sheaves of vector spaces. In \cite{bbd},  the \emph{perverse t-structure} was constructed, whose heart is the abelian category of perverse sheaves.  We denote the cohomology objects with respect to the perverse t-structure by $^{p}\cH^i(\cF^\bullet)$.  Each of the $^{p}\cH^i(\cF^\bullet)$ are perverse sheaves.

Given a local system $\bV_{\bC}$ on a dense open subset $U$ of $\B$, there is a unique extension $^\pi\bV_{\bC}$ of $\bV_{\bC}[n]$ from $U$ to $\B$ which is perverse and has no nontrivial subobjects or quotient objects (in the abelian category of perverse sheaves) supported on $D$.  We put $IC(\bV_{\bC})=\!^\pi\bV_{\bC}[-n]$. If $U=\B$, then $D$ is empty and $IC(\bV_{\bC})$ is just $\bV_{\bC}$.  Some authors define $IC(\bV_{\bC})$ to be the perverse sheaf $^\pi\bV_{\bC}$, but for our applications we find it more convenient to shift degrees so that the lowest nonvanishing cohomology sheaf of $IC(\bV_{\bC})$ is in degree 0.  Both $^\pi\bV_\bC$ and $IC(\bV_\bC)$ can be  represented by complexes of sheaves of vector spaces with constructible cohomology, well-defined up to quasi-isomorphism.

If $X$ is a smooth variety, then $IC(\bR)\simeq\bR$.  We denote this perverse sheaf by $IC(X)$.  More generally, if $X$ is an irreducible variety, let $X^{\rm sm}\subset X$ denote the smooth locus of $X$, a dense open subset of $X$.  Consider the trivial local system $\bR_{X^{\rm sm}}$ on $X^{\rm sm}$.  We then put $IC(X)=IC(\bR_{X^{\rm sm}})$.  Furthermore, if $U\subset X$ is a Zariski open subset of a smooth variety $X$ and $\bV$ is a local system on $U$, then $IC(\bV|_U)\simeq\bV$.  In other words, any local system on a smooth variety can be recovered from its restriction to any Zariski open subset $U\subset X$ by the IC construction.

\subsection{Known results relating $L^2$ and intersection cohomology}\label{subsec:known}

\smallskip
Let $L_{(2)}^p(\bV)$ be the sheaf on $\B$ of measurable $p$-forms valued in $\bV$ which are locally $L^2$ with locally $L^2$ derivatives.  Then the $L_{(2)}^p(\bV)$ fit together to form a complex $L_{(2)}^\bullet(\bV)$.
The following two propositions are among the main results of \cite{cks2,KK}.

\begin{prop}\label{prop:l2ic}
The complex $L^*_{(2)}(\bV_{\bC})$ is quasi-isomorphic to $IC(\bV_{\bC})$.  
\end{prop}

\begin{cor}\label{cor:l2ih}
$H^*(\Gamma(\B,L^*_{(2)}(\bV_{\bC}))) \simeq \bH^*(\B,IC(\bV_{\bC}))$.
\end{cor}
In \cite{cks2,KK} and elsewhere, the hypercohomology groups $\bH^*(\B,IC(\bV_{\bC}))$ are referred to as the intersection cohomology groups $IH^*(\B,\bV_{\bC})$.

It is also shown in \cite{cks2,KK} that there is a good harmonic theory.  In particular
\begin{prop}\label{prop:l2harm}
The cohomology classes in $H^*(\Gamma(\B,L^*_{(2)}(\bV_{\bC})))$ are represented by harmonic forms which are globally $L^2$.
\end{prop}
As an immediate corollary, we have

\begin{cor}\label{cor:harmonic}
\begin{sloppypar}
The vector space of $L^2$ harmonic forms valued in $\bV_{\bC}$ is isomorphic to $\bH^*(\B,IC(\bV_{\bC}))$.
\end{sloppypar}
\end{cor}

While this corollary is all that we really need, it may be helpful at this point to review how to calculate $IC(\bV_{\bC})$, following \cite{cks2}. 

The calculation is local, so we can compute in a neighborhood $W\simeq\Delta^n$ as above with $W\cap U\simeq(\Delta^*)^k\times \Delta^{n-k}$.  A local system $\bV_{\bC}$ on $(\Delta^*)^k\times \Delta^{n-k}$ is equivalent to the data of a finite-dimensional vector space $V$ with monodromies $T_i$ about $z_i=0$ for $i=1,\ldots k$.  In particular, the factor of $\Delta^{n-k}$ is irrelevant and we can reduce to the case $k=n$, $W\simeq \Delta^n$ and $W\cap U \simeq (\Delta^*)^n$.  Near any point $p\in W$ we can find a smaller neighborhood $W_p\subset W$ of $p$ with $W_p\cap U\simeq (\Delta^*)^{k}\times \Delta^{n-k}$ for some $k$ and so we can similarly reduce to the situation $W\cap U\simeq (\Delta^*)^{k}$.  By induction on $n$, we only need to describe the stalk $IC(\bV_\bC)_0$ of $IC(\bV_\bC)$ at the origin.  Here $IC(\bV_\bC)_0$ is a complex of complex vector spaces with finite-dimensional cohomologies which we next describe.

There is a branched cover of $\pi:(\tilde{\B},\tilde{D})\to(\B,D)$ such that the monodromies of $\pi^*(\bV_\bC)$ around the components of $\tilde{D}$ are unipotent \cite{griffiths}.  Both the $L^2$ cohomology and the intersection cohomology on $\B$ are the invariants under covering transformations of the $L^2$ cohomology and the intersection cohomology on $\tilde{\B}$, respectively.  So we may assume that the monodromies are all unipotent.  In our applications to F-theory, the monodromies will be unipotent anyway.

For the local calculation, we consider $\bV_\bC$ to be a local system on $(\Delta^*)^n$ with unipotent monodromy $T_i$ about $z_i=0$.  We think of $T_i$ as an automorphism of a typical fiber $V$.  Then the $N_i=\log(T_i)$ are commuting nilpotent endomorphisms of $V$.  For $0\le p\le n$ we put, following \cite{cks2}
\begin{equation}\label{eq:bp}
    B^p(N_1,\ldots,N_n;\bV_\bC)=\bigoplus_{1\le j_1<\cdots<j_p\le n}N_{j_1}N_{j_2}\cdots N_{j_p}V,
\end{equation}
where conventionally we understand $B^0(N_1,\ldots,N_n;\bV_\bC)=V$.  These vector spaces form the terms of a complex $B^\bullet(N_1,\ldots,N_n;\bV_\bC)$ whose differentials are given on the summands of $B^{p-1}(N_1,\ldots,N_n;\bV_\bC)$ by
\begin{equation}\label{eq:cksdiff}
\left(-1\right)^{s-1}N_{j_s}:N_{j_1}N_{j_2}\cdots\hat{N_{j_s}}\cdots N_{j_p}V \rightarrow N_{j_1}N_{j_2}\cdots N_{j_p}V.
\end{equation}

\begin{prop}\label{prop:ic}
$IC(\bV_\bC)_0$ is represented by the complex $B^\bullet(N_1,\ldots,N_n;\bV_\bC)$.
\end{prop}

We have an immediate corollary.

\begin{cor}\label{cor:icjstar}
$\cH^0(IC(\bV_\bC))\simeq j_*\bV_\bC$.
\end{cor}

To prove this corollary, Proposition~\ref{prop:ic} shows that
\begin{equation}
    \left(\cH^0(IC(\bV_\bC))\right))_0=\cap_{i=1}^n \ker(N_i),
\end{equation}
the space of monodromy invariants.  As $(j_*\bV)_0$ is also the space of monodromy invariants, the corollary follows immediately after considering the reasoning which allowed us to reduce to the stalk at the origin earlier. 

Alternatively, Corollary~\ref{cor:icjstar} can be shown directly and in greater generality using the construction of $IC(\bV)$ from a stratification of $\B$ following \cite[Section~4.2]{dmitm}.

\subsection{Examples}\label{subsec:mathexamples}

\noindent
{\bf Example.} $n=1$, the situation relevant for 8-dimensional F-theory associated to an elliptically fibered K3.  We have a local system $\bV_\bC$ on $\Delta^*$ which is equivalent to the data of a finite-dimensional vector space $V$ with an automorphism $T$.  In our situation, $T$ is unipotent, and $N=\log T$.  Then $B^\bullet(N;\bV_\bC)$ is just
\begin{equation}\label{eq:b1}
V \stackrel{N}{\rightarrow} NV,
\end{equation}
which is quasi-isomorphic to the vector space $\ker(N)$ thought of as a complex in degree 0.
Our stratification is given by $\Delta^*$ and $\{0\}$.  On $\Delta^*$, $IC(\bV)$ is just $\bV$ itself, and on the origin $IC(\bV)_0$ is $\ker(N)$ as just remarked.  It follows that $IC(\bV)$ is represented by a sheaf of vector spaces rather than a complex of sheaves of vector spaces.  Taking Corollary~\ref{cor:icjstar} into account, it follows that 
\begin{equation}
   IC(\bV_\bC)\simeq j_*\bV_\bC. 
\end{equation}
This is the situation in \cite{Zucker}.  In particular, we see Propositions~\ref{prop:l2ic} and \ref{prop:l2harm}, as well as Corollary~\ref{cor:harmonic}, specialize in dimension~1 to the corresponding results in \cite{Zucker}, where $j_*\bV_\bC$ replaces $IC(\bV_\bC)$.

\smallskip\noindent
{\bf Example.} $n=2$.  We have a local system $\bV_\bC$ on $(\Delta^*)^2$ which is equivalent to the data of a finite-dimensional vector space $V$ with automorphisms $T_1$ and $T_2$.  In our situation, the $T_i$ are unipotent, and $N_i=\log T_i$.  Our strata are $(\Delta^2)^2$,
$\Delta^*\times\{0\}$, $\{0\}\times \Delta^*$, and $\{0\}$.  

On $(\Delta^*)^2$, $IC(\bV_\bC)$ is just $\bV_\bC$ itself.  Near $\Delta^*\times\{0\}$, we describe $\bV_\bC$ as a vector space $V$ with monodromy $T_2$ around $z_2=0$ we have already calculated $IC(\bV_\bC)$ as $j_*\bV_\bC$ here by the previous example.  The calculation of $IC(\bV_\bC)$ near $\{0\}\times\Delta^*$ is the same, and we conclude that
\begin{equation}\label{eq:n2ic0}
    IC(\bV_\bC)|_{\Delta^2-0}\simeq j_*\bV_\bC|_{\Delta^2-0}.
\end{equation}
It remains to see how $IC(\bV)$ extends over $0$, and this is given by the complex $B^\bullet(N_1,N_2;V)$, which reads
\begin{equation}\label{eq:b2}
V \stackrel{\phi:=(N_1,N_2)^\mathrm{t}}{\longrightarrow} N_1V\oplus N_2V\stackrel{\psi:=(N_2,-N_1)}{\longrightarrow} N_1N_2V.
\end{equation}

Since $\psi$ is surjective, we see that (\ref{eq:b2}) can be replaced by the quasi-isomorphic two-term complex
\begin{equation}\label{eq:ic2term}
    V \stackrel{\tilde\phi}{\rightarrow}\ker \psi,
\end{equation}
where $\tilde\phi$ is the same map as $\phi$ but with its target restricted.  In situations where $\phi$ is not surjective (and we will see an example presently), we conclude that $IC(\bV)$ is not a sheaf but can be represented a two-term complex $\rho:F^0\to F^1$ 
\begin{equation}
    F^\bullet: F^0\stackrel{\rho}{\rightarrow} F^1
\end{equation}
of sheaves of vector spaces on $\Delta^2$.  We have that 
\begin{equation}
    \cH^0(F^\bullet)\simeq \cH^0(IC(\bV_\bC)) \simeq j_*\bV_\bC.
\end{equation}
Since $IC(\bV_\bC)$ restricts to a sheaf on $\Delta^2-0$ rather than a complex by (\ref{eq:n2ic0}), we have $\cH^1(IC(\bV)))|_{\Delta^2-0}=0$.  Thus $\cH^1(F^\bullet)\simeq \cH^1(IC(\bV_\bC))$ is a skycraper sheaf supported at the origin, where it has stalk isomorphic to the cokernel of $\tilde\phi$.  

We illustrate with an example from 6-dimensional F-theory, where we have a normal crossings intersection of components of $\Delta$ parametrizing fibers of Kodaira type $I_{n_1}$ and $I_{n_2}$, respectively.  Letting
\begin{equation}
    N=\left(
    \begin{array}{cc}
    0&1\\
    0&0
    \end{array}
    \right),
\end{equation}
we have $N_1=n_1N$ and $N_2=n_2N$, and $N_1V$ and $N_2V$ are the same 1-dimensional subspace of $V$.  Referring to (\ref{eq:b2}) we see that $\psi=0$ so that $\ker\psi=N_1V\oplus N_2V$, and finally $\tilde\phi$ has 1-dimensional image.  Thus $\cH^1(IC(\bV_\bC))$ is a skyscraper sheaf at the origin, where it has 1-dimensional stalk.

In particular, this example shows that if $\dim \B > 1$, then $IC(\bV)$ need not be equal to $j_*\bV$.

\section{Physics of F-theory}\label{sec:physics}

F-theory is a geometric approach to studying the physics of string
compactifications.  For the purposes of this paper we can think of
F-theory simply as a procedure for reducing classical type IIB
supergravity from ten dimensions on a compact K\"ahler space in a
regime where string theory is nonperturbative.

In this section we give a brief introduction to the physics of the
systems of interest, describe a conjecture about the structure of the
physical fields, and summarize the relevant mathematics of the
remainder of the paper.

\subsection{Type IIB supergravity: fields and action}

As reviewed more briefly in the Introduction,
Type IIB supergravity is defined by a set of fields ${\cal F}$ on a
ten dimensional manifold with Lorentzian signature, and an action
functional $S({\cal F})$.  The fields ${\cal F}$ include in particular
the space-time metric $g_{MN}$,
the complex scalar {\it axiodilaton} $\tau=\chi+ie^{-\phi}\equiv
\tau_1 + i \tau_2$, and a
pair of antisymmetric two-form fields $B_{\mu \nu}, C_{\mu \nu}$.
The quantum supergravity theory is defined by a path integral over all
field configurations $\int [{\cal DF}] e^{iS({\cal F})}$; the
classical theory is described by the saddle points where the action
functional is stationary $\delta S/\delta f = 0 \forall f \in{\cal
  F}$.

The part of the action (in ``Einstein frame'') that controls the
metric, the
$B$
fields and the axiodilaton is given by
\begin{equation}
S_B = \frac{1}{ \kappa_{10}^2}  \int d^{10} x \sqrt{-g}
\left( R -
{\cal M}_{IJ}F_3^I\cdot F_3^J 
 -\frac{\partial_M \tau \partial^M \bar{\tau}}{ 2 \tau_2^2} 
\right)
\,,
\label{eq:action}
\end{equation}
where 
\begin{equation}
\left(\begin{array}{c}
F_3^1 \\F_3^2
\end{array}
\right)= \left(\begin{array}{c}
dB\\dC
\end{array} \right) \,,
\label{eq:?}
\end{equation}
and
\begin{equation}
{\cal M}_{IJ} =\frac{1}{\tau_2}  \left(\begin{array}{cc}
| \tau |^2 & -\tau_1\\
-\tau_1 & 1
\end{array} \right) \,.
\label{eq:m-physics}
\end{equation}
There is also a four-form field with additional terms in the action
but for our purposes we can assume that this field vanishes and these
terms do not contribute.
We will also assume that $B, C = 0$ in the vacuum solutions of
interest but we will consider the fluctuations of the two-form fields
around this vacuum, which constitute the physical degrees of freedom
of interest.

\subsection{Dimensional reduction of IIB to 8D and 6D}

We begin by dimensionally reducing the type IIB supergravity
theory to 8 dimensions by
compactifying on ${\cal B} =\bP^1$.  
We assume that the fields can vary as functions of a coordinate $z$ on
${\cal B}$, and are independent of the remaining space-time coordinates.
The equations of motion for the axiodilaton $\tau$ are independent of
the metric and dictate that $\tau$ is locally a holomorphic function
of $z$.
When the axiodilaton $\tau$ does not depend
upon the position in $\bP^1$,
the equations of motion for the metric
$g_{MN}$ are that locally the metric is Ricci flat.  There are,
however, extended objects known as $(p, q)$ 7-branes that carry
$\pi/6$ units of curvature, and which source the axiodilaton; $\tau$
has a corresponding monodromy
around any such 7-brane.
Mathematically, these can be thought of as poles/singularities in the
axiodilaton $\tau$.  
Clearly on $\bP^1$ there must be 24 7-branes for a total curvature of
$4 \pi$, and the product of the monodromies must be trivial in an
appropriate coordinate system.
The axiodilaton is characterized by 
a Weierstrass model of the form (\ref{eq:weierstrass}), where $f, g$
are degree 8, 12 polynomials in $z$; this gives
an elliptic fibration over $\bP^1$, and the resulting axiodilaton configuration
 can be described in terms of the Weierstrass coefficients $f, g$
 through
\begin{equation}
j (\tau (z)) = 1728 \times \frac{4f^3}{\Delta} , \;\;\;\;\;
\Delta = 4f^3 + 27g^2 \,.
\end{equation}
The metric obeys the Einstein
equations with the variations of $\tau$ providing additional local
source terms from the action (\ref{eq:action})
\cite{Greene:1989ya, Gibbons:1995vg}.  Writing an Ansatz for
the metric in
the form
\begin{equation}
ds^2 = \sum_{i = 0}^{7}  dx_i^2 + \Omega (z, \bar{z}) dz\ d\bar{z} \,,
\label{eq:metric}
\end{equation}
the equation of motion for the metric is
\begin{equation}
\partial \bar{\partial} \tau_2 = 2\partial \bar{\partial} \Omega \,.
\end{equation}
It follows that $\Omega = \tau_2 e^{F + \bar{F}}$ where $F (z)$ is
locally a holomorphic function. In the vicinity of a seven-brane
localized at e.g. $z = 0$ the resulting metric is conformally
equivalent to the flat metric $dz d\bar{z}$ in a neighborhood of $z =
0$ on the punctured disk.
In \cite{Einhorn:2000ct}, a further analysis of the structure of seven-brane metrics
is given, showing that the metric from supergravity
becomes singular already at a finite but small distance away from the
seven-brane where the supergravity approximation breaks down.
Nonetheless, the normalizability conditions on one-form fields in this
metric are the same as for (\ref{eq:metric}), at least for 8D
compactifications, since the metrics are conformally equivalent and
the factors of the metric cancel in computing the norm.

The story is similar for compactification to 6D; the compactification
space $\B_2$ is a compact K\"ahler surface, and the 7-branes are
described by a divisor in the class $-12K_{B_2}$.
In this case, however, the seven-brane locus generally contains
codimension two singularities that complicate the structure of the metric.

We are interested in fluctuations in the $B, C$ fields around their
vanishing values in the vacuum.  In particular, when we can write
$B = \phi \wedge A$ with $\phi$ a one-form in the compact space and
$A$ a one-form in the non-compact directions, then when $d \phi = 0$
the action reduces to the Maxwell action $dA^2$ in the dimensionally
reduced theory.  When $d \phi \neq 0$, then there is a mass term
$A^2$.  We are interested in light fields in the reduced theory so
focus on reduction with $d \phi = 0$.  Thus, each closed one-form on
the compactification space $\bP^1$ or $\B_2$ gives rise to a Maxwell
(U(1)) gauge field in the 8D or 6D reduced theory.
Similarly, each closed two-form gives a scalar field in the 8D or 6D
theory.

\subsection{A conjecture regarding normalization}\label{subsec:normalization}

String theory is an approach to giving a quantum description of
supergravity.  Type IIB string theory is a perturbative theory defined
in the regime where the axiodilaton $\tau$ is constant in the full
10-dimensional space-time, and the string coupling $e^{-\phi}$ is small.
From the point of view of string theory, the 7-branes of the
type IIB theory
 are nonperturbative dynamical objects.  In the
presence of a D7-brane, which gives monodromy $\tau \rightarrow \tau +
1$ but leaves $e^{-\phi}$ invariant, the dynamics of string theory is
described by a combination of closed strings without endpoints that
propagate in the ``bulk'' space-time away from the D7-brane and open
strings that end on the D7-brane.  From the point of view of string
theory, the closed strings give fields like the space-time graviton
while the open strings give gauge fields that propagate on the brane.
A single isolated brane carries a U(1) gauge factor, while
multiple D7-branes that are coincident give a nonabelian gauge field
$U(N)$.
A similar story should hold for the dynamics of other types of
7-branes, which are related by SL(2,Z) duality to the D7-brane, though
there is no perturbative description of these dynamics.

The notion of holography, associated with the AdS/CFT correspondence,
asserts that there is a duality between the degrees of freedom
described by open strings as the gauge field on a brane and the
dynamics of the supergravity fields in the ``near horizon'' region of
geometry right next to the brane.  In this paper, inspired by the 8D
analysis of Douglas, Park, and Schnell, we speculate that from the
type IIB point of view, the relevant dynamical fields of the
compactification space ${\cal B} \setminus \Delta$ can be separated
into two sets of degrees of freedom: those that are normalizable 
and
those that are not.  We conjecture that this decomposition matches
with the separation of bulk and brane degrees of freedom expected from
string theory.  In particular, stated in  somewhat imprecise  physical terms,

\begin{conj}
Closed one-forms in ${\cal B}\setminus \Delta$ that are normalizable
correspond to abelian gauge fields in the reduced 8D or 6D theory,
while non-normalizable one-forms correspond to gauge fields living on
7-branes localized on the discriminant locus $\Delta$.  Similarly,
normalizable two-forms in ${\cal B}\setminus \Delta$ correspond to
scalars in the reduced theory coming from bulk supergravity dynamics
while non-normalizable two-forms correspond to scalars localized on
the 7-branes.
\label{c:physics}
\end{conj}

In general, the 
abelian gauge degrees of freedom in $B, C$ associated with
seven-brane dynamics  (i.e. by the conjecture those that are
non-normalizable) are removed by gauge redundancy through the
Cremmer-Scherk mechanism as in the 8D case; for example,  if the
discriminant locus is purely of $I_1$ type there are no nonabelian
factors and the only abelian gauge factors are those associated with
the bulk dynamics, while if there is an $I_N$ locus the associated
$\U(N)$ gauge field on the brane is reduced by Cremmer-Scherk to an
$\SU(N)$ gauge field.
In general we expect that the (conjectured normalizable) scalars in
the bulk will correspond to uncharged hypermultiplets, while the
(conjectured non-normalizable) scalars localized on $\Delta$ will
generally correspond to scalars charged under the gauge group.  In
this paper we do not address questions related to the charged scalars,
only the  gauge fields and uncharged scalars, leaving further
treatment of the charged scalars to further work.

Translating this conjecture into the language of mathematics, 
in analogy to de Rham cohomology on a compact manifold,
we
basically expect that the complete set of closed
$p$-form fields on ${\cal B} \setminus \Delta$ (for $p = 1, 2$)
modulo exact fields
will live in some cohomology theory $H^p_{\rm all} ({\cal B} \setminus
\Delta)$, while the normalizable fields will live in another
cohomology theory  $H^p_{\rm norm} ({\cal B} \setminus
\Delta)$, and the non-normalizable fields associated with the boundary
will live in  $H^{q}_{\rm  boundary} (\Delta)$ where we are agnostic
for the moment about the value of $q$.  We are then
looking for a mathematical formulation of a set of cohomology theories
that will give us an exact sequence of the form
\begin{equation}
0 \rightarrow
H^p_{\rm norm} ({\cal B} \setminus \Delta)\rightarrow
H^p_{\rm all} ({\cal B} \setminus \Delta)\rightarrow
H^{q}_{\rm  boundary} (\Delta)\rightarrow
0
\label{eq:bulkboundary}
\end{equation}
Such an exact sequence amounts to a decomposition of the full set of
fields into a normalizable set of fields, which is a subset of the
full set, and the set of boundary fields, which are the quotient of
the full set by the normalizable set.  

After describing F-theory and some specific examples in a little more
detail, the remaining sections of this paper attempt to formulate a
precise mathematical theory that will realize this decomposition of
the physical degrees of freedom between boundary and bulk, potentially
giving insight not only into F-theory compactifications but also into
aspects of the nature of holography.  

It may be helpful here to briefly summarize some of the main
mathematical
points in
the remainder of the paper and how they connect to this schematic
physics picture.
A principal focus is the set of  $p$-forms with
a normalizable representative with respect to the physical metric
coming from dimensional reduction of IIB supergravity.
In 8D compactifications, as argued in \cite{dps} each normalizable
$1$-form field has a harmonic representative with minimal norm.
For the Poincar\'{e} metric (\ref{eq:asympmetric}),
the cohomology of these fields is given by $H^1(\bP^1,j_*\bV)$, which
is isomorphic to     $H^1(\bP^1,IC (\bV))$.
The primary conjecture (Conjecture \ref{c:main})
in the following section asserts that more
generally, for the physical cases of interest, what we have called
$H^p_{\rm norm} ({\cal B} \setminus \Delta)$ here is given by
$H^p({\cal B},IC (\bV))$.  
The full set of ($\bC^\infty$) $p$-form cohomology classes on
$U ={\cal B} \setminus \Delta$ is given by $H^p(U,\bV)$.
The expectation then is that $H^p({\cal B},IC (\bV))$ describes the
set of physically relevant ``bulk'' fields, that
$H^p(U,\bV)$ describes the full set of fields, and that
there is an injective map from the former to the latter, with the
quotient localized on the 7-branes.

In 8D this injective map exists and is part of the exact sequence
(\ref{eq:holseq}).  While the class of metrics such as
(\ref{eq:asympmetric}) for which $H^1(\bP^1,j_*\bV)$ gives
normalizable forms is different from the specific form of the metric
(\ref{eq:metric}), for 8D theories it shares the conformally flat form
of the metric.  Furthermore, the inner product defined by
(\ref{eq:m-physics}) matches precisely with that used in
\cite{Zucker}.  Thus, with this inner product an $\SL(2,\bZ)$ doublet
of harmonic forms is normalizable with respect to the metric
(\ref{eq:asympmetric}) if and only if it is normalizable with respect
to (\ref{eq:metric}).  Since such  an $L^2$ harmonic form is
$\bC^\infty$, the injection from
 $H^1(\bP^1,j_*\bV)\rightarrow H^1(U,\bV)$ allows us to identify the
physically relevant doublet one-form fields from the bulk
with a subset of the full set of forms.  This analysis, which is
essentially that of \cite{dps}, shows that the decomposition of forms
is correctly described by Conjecture \ref{c:main} for the 8D case.

For 6D theories there are a number of issues that would need to be
resolved to prove Conjecture \ref{c:main}.  In particular, the
structure of the metric for which  $H^p({\cal B},IC (\bV))$ gives
normalizable forms does not necessarily match that coming from
physics, and the normalizability conditions are not obviously
equivalent.  
In general, the discriminant locus $\Delta$ contains not only
transverse intersections but also cusp singularities, for which the
corresponding mathematical analysis is not known.
Furthermore, for higher-dimensional $\B$ we do not have a proof that
the map from 
 $H^p({\cal B},IC (\bV))$ to  $ H^1(U,\bV)$ is injective, although we
discuss some aspects of this possibility in \S(\ref{subsec:supports}).
We leave a more detailed analysis of these issues for
further work.
We do, however, compute the dimension of $H^p({\cal B},IC
(\bV))$ in a number of examples enumerated in Section
\ref{sec:examples}, and find that the resulting numbers are compatible
with the assertions of Conjectures \ref{c:physics} and \ref{c:main}.

\subsection{F-theory}

F-theory \cite{Vafa-F-theory, Morrison-Vafa-I, Morrison-Vafa-II} is a
nonperturbative approach to string theory that describes many features
of type IIB compactifications on general ${\cal B}$ in the presence of
7-branes.  F-theory is usually defined as a limit of the
11-dimensional M-theory version of string theory, and the definition
of F-theory is as yet mathematically incomplete.  Part of the purpose
of this paper is to build our understanding of F-theory directly as a
method for understanding the nonperturbative physics of type IIB
compactification.

F-theory at some level is a dictionary between the geometry of a $d$-dimensional
elliptic Calabi-Yau manifold and the physics of the resulting
compactification to $12-2d$ space-time dimensions.  The elliptic
Calabi-Yau manifold $X$ is described by the axiodilaton $\tau$ given by
 a Weierstrass model over the complex base ${\cal B}$.
Singularities in the elliptic fibration give rise to physical features
in the reduced theory.  Codimension one (in the base ${\cal B}$)
singularities of the elliptic fibration  give rise to nonabelian gauge
fields in the reduced theory where the gauge algebra is captured by
the Dynkin diagram of the associated Kodaira singularity type.
Codimension two singularities give rise to matter fields.
The number of connected abelian (U(1)) factors in the gauge group is
given by the rank of the Mordell-Weil group of sections of the
elliptic fibration $X$.
While this statement about the abelian part of the gauge group is
usually argued for from the physics of M-theory, here we give (in \S\ref{sec:MW}) a direct
demonstration  of this relationship in general dimension by showing
that a section of $X$ gives a doublet of one-forms in ${\cal B}$ that
give the appropriate U(1) factor in the reduced theory.

For compactifications on a complex surface ${\cal B}$ (the case $d =
2$), the resulting relations between geometry and physics nicely match
with the expectation from gravitational and gauge anomaly cancellation
in six-dimensional supergravity theories.
Some features of the physics and the dictionary to geometry in this case include
the following:
\begin{eqnarray}
T & = & h^{1, 1} ({\cal B}) -1\\
h^{1, 1}(X) & = & h^{1, 1} ({\cal B}) + 1 + {\rm rk}\ G\label{eq:stw-physics}\\
H_{\rm unch} & = & 1 +h^{2,1} (X)  \\
H-V &= & 273-29T\,, 
\end{eqnarray}
where $T$ is the number of tensor multiplets of the 6D
supergravity theory, which contain
antisymmetric 2-tensor fields, $G$ is the (abelian
and nonabelian) gauge group of the 6D theory, $V= {\rm dim}\ G$ is the
number of vector multiplets, and $H, H_{\rm unch}$ are the numbers of
(all, uncharged respectively) scalar hypermultiplets in the 6D theory.
The last of these relations expresses the gravitational anomaly
consistency condition of 6D supergravity.  Note that the uncharged
scalar hypermultiplet fields have four real degrees of freedom; two
come from the axiodilaton $\tau$ and characterize the complex
structure moduli of the elliptic Calabi-Yau, while the other two real
components come from $B$ fields and are the focus of the analysis
here.

\subsection{Examples}\label{sec:examples}

We conclude this section with a few simple examples of 6D F-theory
compactifications that we consider further from the mathematics
perspective in the next section.
\vspace*{0.05in}

\noindent {\bf Base ${\cal B} =\bP^1$}

In this case $X$ is an elliptic K3 surface.  This is the case treated
by DPS in \cite{dps}.  The rank of the gauge group is 20, and can
either come from a purely abelian group $\U(1)^{20}$ or can include
nonabelian factors from Kodaira singularities in the elliptic
fibration.  The primary result of the DPS paper was to derive the rank
of the gauge group as 20 in the purely abelian case directly from the
IIB perspective.
\vspace*{0.05in}

\noindent {\bf Base ${\cal B} =\bP^2$}

In this case $X$ is the generic elliptic fibration over the projective
surface $\bP^2$, which has Hodge numbers $(h^{1, 1},h^{2, 1}) = (2,
272)$.  In this case we thus expect that there will be  273 uncharged
scalars in this 6D theory.  
\vspace*{0.05in}

\noindent {\bf Gauge group $\SU(2)$ over a line in ${\cal B} =\bP^2$}

In this case, $X$ is a Weierstrass model over ${\cal B} =\bP^2$, which
has been tuned to have a Kodaira $I_2$ singularity over a line in the
projective base.  This model has 22 codimension two points in the base
where the $I_2$ locus intersects the residual $I_1$ locus, each
corresponding to a hypermultiplet of scalar fields in the fundamental
representation of the $\SU(2)$ gauge group.  The gravitational anomaly
condition $H-V= 273-29T$ (with $T = 0, V= 3$)
then tells us that $H_{\rm unch} =232$.

\section{The main conjecture and application to F-theory}
\label{sec:main-conjecture}

We return to our main interest, elliptically fibered Calabi-Yau manifolds.

\subsection{Statement of the conjecture}

Let $\pi:X\to \B$ be an elliptically fibered Calabi-Yau manifold with discriminant divisor $\Delta\subset \B$.  We put $U=\B-\Delta$,  $\bV=R(\pi_U)_*\bR$, and $\bV_\bC=\bV\otimes_\bR\bC\simeq R(\pi_U)_*\bC$.

Suppose first that $S\to \B=\bP^1$ is an elliptically fibered K3
surface with 24 nodal fibers (Kodaira type $I_1$).  This is the
situation considered in \cite{dps}.  Since the monodromy is unipotent,
\cite{Zucker} applies (or we can apply \cite{cks2,KK} with $n=1$,
noting that $IC(\bV_\bC)\simeq j_*(\bV_\bC$) by the $n=1$ example in
Section~\ref{subsec:mathexamples}.  We conclude that the space of harmonic
1-forms valued in $\bV_\bC$ is
isomorphic to $H^1(\bP^1,j_*\bV_\bC)$. 
As mentioned above, this space of harmonic forms is the same for the
physical metric and the Poincar\'{e} metric.

Via this isomorphism, complex conjugation acts compatibly on harmonic 1-forms valued in $\bV_\bC$ and on $H^1(\bP^1,j_*\bV_\bC)$.  We conclude that

\medskip\noindent
The space of real harmonic forms valued in $\bV$ is isomorphic to    $H^1(\bP^1,j_*\bV)$.

\medskip\noindent This gives a mathematical proof that the vector space of gauge fields
arising from dimensionally reducing the B-fields is isomorphic to
$H^1(\bP^1,j_*\bV)$.  Shortly, we will show how to use the
Decomposition Theorem of \cite{bbd} to compute that $\dim
H^1(\bP^1,j_*\bV)=20$, as was anticipated by the M-theory description.

For $n>1$, the discriminant divisor $\Delta$ is typically not normal crossings due to the presence of cusps. Also, there is no known simple form for the asymptotics of the metric in complete generality. Nevertheless, we have a conjectural way forward.

\begin{conj}
\label{c:main}
Consider an F-theory compactification associated to an elliptically
fibered Calabi-Yau $\pi:X\to \B$ in any dimension.  Let $\bV$ be
the local system $R^1(\pi_U)_*\bR$ on $U=\B-\Delta$.  Then the vector space of $\bV$-valued $i$-forms on $U$, harmonic and $L^2$ with respect to the K\"ahler metric, is isomorphic to $\bH^i(\B,IC(\bV))$.
\end{conj}

Assuming this conjecture we get immediately

\begin{cor}\label{cor:physicsofB}
\ \begin{enumerate}
    \item The space of abelian
gauge fields obtained by dimensionally reducing normalizable B-fields is parametrized by $\bH^1(X,IC(\bV))$.
\item The space of uncharged
scalar fields obtained by dimensionally reducing normalizable B-fields is parametrized by $\bH^2(X,IC(\bV))$.
\end{enumerate}
\end{cor}

In the rest of this paper, we will check our conjecture by showing
that it is consistent with results of physics, including the dual
M-theory description, assuming that Conjecture~\ref{c:physics}
correctly captures the mathematical conditions on physical
fields.  
Specifically, we check consistency of the first
statement of Corollary~\ref{cor:physicsofB} with physics for all $X$
of dimension 3 which admit a crepant resolution and whose fibers over
generic points of the components of $\Delta$ have Kodaira type
$I_n,\ I_n^*,\ II^*,\ III^*$, or $IV^*$ and such that
$H^{2,0}(X)=H^{2,0}(B)=0$.  We also check consistency of the second
statement of Corollary~\ref{cor:physicsofB} for all smooth $X$ of
dimension 3 over a regular base $B$ with $I_1$ fibers over a generic
point of $\Delta$.

\subsection{The Decomposition Theorem}

We recall the Decomposition Theorem of \cite{bbd}, beginning with some background.

\smallskip
Any object $F^\bullet$ of any derived category has the same cohomology objects $\cH^i(\cF^\bullet)$ as does $\oplus_i\cH^i(\cF^\bullet)[-i]$, but there is no reason to expect these be the same objects of the derived category.  However, a theorem of Deligne says that we  have such an isomorphism in an important situation.
\begin{prop}\label{eq:delignethm}
Suppose that $f:X\to Y$ is a proper smooth morphism of quasi-projective varieties of relative dimension $m$.  Then
\begin{equation*}
    Rf_*\bR_X\simeq \bigoplus_{i=0}^m R^if_*R_X[-i].
\end{equation*}
\end{prop}
The sheaves $R^if_*R_X$ appearing in (\ref{eq:delignethm}) are the local systems on $Y$ whose stalks at $y\in Y$ are the cohomologies $H^i(X_y,\bR)$ of the fibers $X_y$ of $f$ over $y$.

If in addition the fibers of $f$ are irreducible, then $H^0(X_y,\bR)$ and $H^{2m}(X_y,\bR)$ are canonically trivial, so that $R^0f_*\bR_X\simeq R^{2m}f_*\bR_X\simeq \bR_Y$.

\smallskip
We now recall the Decomposition Theorem, which says that
Proposition~\ref{eq:delignethm} extends to non-smooth $f$ if we use
the perverse t-structure instead of the standard t-structure.  For an
object $F^\bullet$ of the constructible derived category of $X$ and a
map
\begin{equation}
    ^pR^if_*F^\bullet ={}^p\cH^i\left(Rf_*F^\bullet\right),
\end{equation}
where $^p\cH^*$ denotes perverse cohomology, the cohomology with respect to the perverse t-structure. 

\begin{thm}
Let $X$ be an irreducible projective variety and let $f:X\to Y$ be a proper morphism.  Then 
\begin{equation*}
    Rf_*IC(X)=\bigoplus_i  {}^pR^if_*IC(X)[-i],
    \end{equation*}
 Furthermore, each perverse sheaf\ \ $^pR^if_*IC(X)$ appearing in the Decomposition Theorem is the direct sum of the perverse shifts $^\pi\bL$ of the complexes $IC(\bL)$ associated to local systems $\bL$ on dense open subsets of closed subvarieties $Z\subset Y$.
\end{thm}
We call the closed subvarieties $Z$ appearing in the Decomposition Theorem the \emph{supports} of $Rf_*IC(X)$.

\subsection{The Decomposition Theorem in 8-dimensional F-theory}

We consider an elliptically fibered K3 surface $\pi:S\to \bP^1$ with 24 nodal fibers and identify the dimensional reductions of the B-fields.  In \cite{dps}, 
it was computed that $\dim H^1(\mathbb{P}^1,j_*\bV)=20$, which is then
identified with the dimension of the space of gauge fields, i.e.\ the
rank of the gauge group.  This is consistent with the usual
understanding of F-theory as a limit of M-theory, which says that the rank of the gauge group is $\dim H^{1,1}(S)=20$.

Rather than sketch the computation of \cite{dps}, instead we compute $H^1(\bP^1,j_*\bV)$ using the  Decomposition Theorem, which is better suited for generalizations.   Since the restriction of $\pi$ to $\pi|_U:\pi^{-1}(U)\to U$ is smooth, Deligne's Theorem applies and 
\begin{equation*}
    (R\pi_*\bR_S)|_{U}\simeq \bR_U\oplus \bV[-1] \oplus \bR_U[-2].
\end{equation*}
We want to extend this from $U$ to all of $\bP^1$. The result is \cite[Example~1.8.4]{dmdpt}
\begin{equation}\label{eq:k3decomp}
R\pi_*\bR_S \simeq \bR_{\bP^1}\oplus j_*\bV[-1]\oplus \bR_{\bP^1}[-2].
\end{equation}

For later use, we note that \cite[Example~1.8.4]{dmdpt} says the following more generally.  Let $f:S \to C$ be a projective map with connected fibers from a smooth
surface $S$ onto a smooth curve $C$ . Let $\Sigma\subset C$ be the finite set of critical values of $f$ and let
$U = C\backslash \Sigma$ be its complement. The map $f$ is a $C^\infty$ fiber bundle over $U$ with typical fiber a compact Riemann surface of some fixed genus $g$. Let $\bV = (R^1f_*\bR)|_U$ be the rank $2g$ local
system on $U$ with stalk the first cohomology of the typical fiber. We then have an isomorphism
\begin{equation}\label{eq:mdccurves}
Rf_*\bR_S\simeq \bR_C\oplus j_*\bV[-1]\oplus T_\Sigma[-2] \oplus \bR_C[-2],    
\end{equation}
where $T_\Sigma$ is a skyscraper sheaf over $\Sigma$ with stalks $T_s\simeq H_2(f^{-1}(s))/\langle |f^{-1}(s)|\rangle$ at $s\in\Sigma$.

Returning to the elliptically fibered K3 $\pi:S\to\bP^1$, we can compute $H^1(\bP^1,j_*\bV)$ by computing the hypercohomology $\bH^2$ of both sides of (\ref{eq:k3decomp}):
\begin{equation}
H^2(S,\bR)\simeq H^2(\bP^1,\bR)\oplus H^1(\bP^1,j_*\bV)\oplus H^0(\bP^1,\bR),
\end{equation}
from which we easily conclude that $\dim H^1(\bP^1,j_*\bV) = 20$.

Similarly, the 8D scalars are described by dimensionally reducing the B-fields on 2-forms valued in $\bV$, which are parametrized by $H^2(\bP^1,j_*\bV)$.  Taking $\bH^3$ of both sides of (\ref{eq:k3decomp}) we obtain
\begin{equation}
H^3(S,\bR)\simeq H^3(\bP^1,\bR)\oplus H^2(\bP^1,j_*\bV)\oplus H^1(\bP^1,\bR),
\end{equation}
from which we conclude that $H^2(\bP^1,j_*\bV)=0$, and no scalars arise from dimensionally reducing the B-fields.

\subsection{The Decomposition Theorem in 6-dimensional F-theory}

We now turn to the case of an elliptically fibered Calabi-Yau threefold $\pi:X\to \B$.  Suppose that we are in the a generic situation where the fibers of $\pi$ are all irreducible elliptic curves, nodal over the smooth locus of $\Delta$ and cuspidal over the cusps of $\Delta$, where $\Delta$ is an irreducible curve having only cusps as singularities.  We assert that in this case, the Decomposition Theorem reads
\begin{equation}\label{eq:cydecomp}
R\pi_*\bR_X \simeq \bR_B\oplus IC(\bV)[-1]\oplus \bR_B[-2].
\end{equation}
To see this, we begin by identifying the possible $Z\subset \B$ appearing as support of $R\pi_*\bR_X$.  Let $\kappa\subset\Delta$ denote the set of cusps of $\Delta$ and put $\Delta^0=\Delta-\kappa$.  Since the fibers of $\pi$ are equisingular over the three strata $\kappa$, $\Delta^0$, and $U$, we see that $R\pi_*\bR_X$ is locally constant on the locally closed subvarieties $\kappa,\ \Delta^0$, and $U$.  So the supports $Z$ can only be among the respective closures $\kappa,\ \Delta$, and $\B$.

Next, $\kappa$ can be ruled out by the Goresky-MacPherson inequality, which says that if $\pi:X\to Y$ is a proper map of algebraic varieties with $X$ smooth and where all fibers have the same dimension $d$, then for each $Z\subset Y$ appearing as the support of a perverse sheaf in the Decomposition Theorem, we have $\mathrm{codim} Z\le d$.  In the case of an elliptic fibration, we have $d=1$, so the codimension 2 locus $\kappa$ is ruled out. A proof of the Goresky-MacPherson inequality can be found in \cite[Th\'eor\`eme 7.3.1]{ngo}.

So the only possible supports are $\Delta$ and $\B$.  Restricting to a generic curve $C\subset \B$ (which in particular does not contain the cusps), we have an elliptic fibration $\pi|_{\pi^{-1}(C)}:{\pi^{-1}(C)}\to C$, and we are back in the situation of \cite[Example~1.8.4]{dmdpt}.  The surface ${\pi^{-1}(C)}$ is smooth, with critical values $\Sigma=C\cap \Delta$.  Since the fibers of $\pi$ over points of $\Sigma$ are irreducible nodal curves, we see that $T_\Sigma=0$ and there are no supports of $(R\pi_*\bR)|_C$ on points, only on all of $C$.  If follows that $R\pi_*\bR$ cannot have any supports on $\Delta$, since the restriction of $R\pi_*\bR$ to $C$ would have supports on $\Sigma$.  So the only supports of $R\pi_*\bR_X$ can be on $\B$ itself, arising from local systems on $U$, and so we only need to identify the local systems on $U$.  However, by Deligne's Theorem we have
\begin{equation}\label{eq:deligne6d}
    \left(R\pi_*\bR_X\right)|_U\simeq \bR_U\oplus\bV|_U[-1]\oplus\bR_U[-2].
\end{equation}
It follows that $R\pi_*\bR_X$ is isomorphic to the direct sum of the IC sheaves of the terms on the right hand side of (\ref{eq:deligne6d}).  This proves (\ref{eq:cydecomp}).

\smallskip
We can now compute the dimensional reduction of the B-fields.  Taking $\bH^2$ of (\ref{eq:cydecomp}) we get
\begin{equation}\label{eq:gengauge}
H^2(X,\bR)\simeq H^2(\B,\bR)\oplus \bH^1(\B,IC(\bV))\oplus H^0(\B,\bR).
\end{equation}

As a special case, suppose that $\B=\bP^2$ and $X\to \bP^2$ is a
generic Weierstrass fibration.  We then have $h^{1,
  1}(X)=2$ 
and $h^2(\bP^2)=1$. So (\ref{eq:gengauge}) implies that $\bH^1(\B,IC(\bV))=0$ and no gauge fields arise from dimensionally reducing the B-fields.

Before returning to a generic elliptically fibered Calabi-Yau threefold, 
we recall the Shioda-Tate-Wazir formula \cite{shioda,tate,wazir},
which holds more generally.  Let $\mathrm{MW}(X)$ be the Mordell-Weil
group of rational sections of $\pi$, let $\Delta_k$ be the irreducible
components of $\Delta$, and suppose that $\pi^{-1}(\Delta_k)$ has
$n_k$ irreducible components.  Then, in accord with (\ref{eq:stw-physics}),
\begin{equation}\label{eq:stw}
    h^{1,1}(X) = 1 + \operatorname{rank}\mathrm{MW}(X) + h^{1,1}(\B) + \sum_k \left(n_k-1\right).
\end{equation}
In our situation, $\Delta$ is irreducible with $\pi^{-1}(\Delta)$ having $n_1=1$ component.  For simplicity, assume that $h^{2,0}(X)=h^{2,0}(\B)=0$, which holds in almost all cases of interest anyway.   We conclude that
\begin{equation}\label{eq:stwgen}
    h^{2}(X) = h^{2}(\B) + \operatorname{rank}\mathrm{MW}(X) + 1.
\end{equation}
Comparing (\ref{eq:stwgen}) and (\ref{eq:gengauge}), we see that the right-hand side of (\ref{eq:stwgen}) is just the sum of the dimensions of the summands in the right-hand side of (\ref{eq:gengauge}).  In particular, our methods imply:
\begin{equation}\label{eq:mwic}
\operatorname{rank}\mathrm{MW}(X)=\dim\bH^1(\B,IC(\bV)).
\end{equation}
In fact, we will see below that in general, consistency of our conjectures with the Shioda-Tate-Wazir formula is equivalent to (\ref{eq:mwic}).

We will return to the role of the Mordell-Weil group in Section~\ref{subsec:themap}, where we will define a natural map $\mathrm{MW}(X) \to \bH^1(\B,IC(\bV))$.  

In the special case $\B=\bP^2$, it is known that $MW(X)=0$ for a generic $X$, consistent with $\bH^1(\bP^2,IC(\bV))=0$.

\medskip
We now turn our attention to the scalar fields.

\smallskip
In the case of general $\B$, taking $\bH^3$ of (\ref{eq:cydecomp}) we get
\begin{equation}
H^3(X,\bR)\simeq H^3(\B,\bR)\oplus \bH^2(\B,IC(\bV))\oplus H^1(\B,\bR).
\end{equation}

Suppose further that $\B$ is regular surface, so that $H^1(\B,\bR)=H^3(\B,\bR)=0$.  Then we conclude that $\bH^2(\B,IC(\bV))\simeq H^3(X,\bR)$, so we have a real $h^3(X)$-dimensional space of scalars arising from dimensionally reducing the B-fields.

This result is entirely consistent with  expectations from F-theory.  The hypermultiplet moduli space has quaternionic dimension $h^{2,1}(X)+1$, and can be identified with a quantum-corrected version of the Calabi-Yau integrable system $\cJ\to\widetilde{\cM}$ \cite{DM}.  Here $\widetilde{\cM}\to \cM$ is a $\bC^*$-bundle over the complex structure moduli space $\cM$ of $X$, with fiber over $X\in \cM$ the space of nonvanishing holomorphic 3-forms $\Omega$ on $X$, and the fiber of $\cJ\to\widetilde{\cM}$ over $(X,\Omega)$ is the intermediate Jacobian $J(X)$. 

We have already observed that the moduli space of F-theory has a
contribution of $\cM$, which has complex dimension $h^{2,1}(X)$.
F-theory has an additional universal complex modulus which we do not
elaborate on here, giving $h^{2,1}(X)+1$ complex moduli.  According to
duality, we are supposed
to
have an additional $h^{2,1}(X)+1$ complex
parameters or equivalently $2(h^{2,1}(X)+1)$ real parameters.
However, since $h^3(X)=2(h^{2,1}(X)+1)$, we see that these parameters
correspond precisely to the scalars obtained by dimensionally reducing
the B-fields.  For example, in the case $\B=\bP^2$, we have
$h^{2,1}(X)=272$ and we get 273 complex scalars from dimensionally
reducing the B-fields.  These observations also confirm the count of
uncharged scalars noted in Section~\ref{sec:examples}. The Hodge
numbers of $X$ are $(2,272)$.

\smallskip
We now turn to the non-generic situation, where $\Delta$ can have
several components, and $X$ not be smooth.  We start with an example
before explaining the general theory.

\smallskip\noindent
{\bf Example.}  $I_2$ along a line in $\bP^2$.

\smallskip
Suppose that we are in the generic situation where we have an $I_2$ fiber along a line $\ell\subset\bP^2$.  We can describe this in Weierstrass form by letting $\ell$ also denote a linear form on $\bP^2$ vanishing along the line, and taking
\[
f=-3\sigma_6^2+\ell f_{11}+\ell^2f_{10}, \qquad g=2\sigma_6^3- \sigma_6\ell f_{11} +\ell^2g_{16},
\]
where the subscripts denote the degree of a polynomial.  We compute
\begin{equation}\label{eq:i2disc}
    \Delta=\ell^2\sigma_6^2\left(-9f_{11}^2+108 g_{16}\sigma_6+108 f_{10}\sigma_6^2
\right)+O(\ell^3).
\end{equation}
The degree~36 discriminant curve factors into the line $\ell$ (with multiplicity 2) and an irreducible curve $\Delta_1$ of degree~34. The curve $\ell$ meets $\Delta_1$ in the 6 points $\sigma_6=0$ (with multiplicity 2) together with the 22 points $-9f_{11}^2+108 g_{16}\sigma_6+108 f_{10}\sigma_6^2=0$.  We denote this set of 22 points by $Z\subset\ell$.

In this situation, the Calabi-Yau threefold is singular along the curve $C\subset X$ given by
\begin{equation}
\ell=0,\ x=\sigma_6,\ y=0,\ z=1.   
\end{equation}
The transverse singularity along $\ell$ is generically an $A_1$ singularity, but the singularity is more degenerate along $Z$.

Since $X$ is singular, we can analyze the B-fields using the Decomposition Theorem applied to $R\pi_*IC(X)$.  However, a simpler computational method becomes evident after resolving $X$.
Blowing up $C$ produces a smooth Calabi-Yau $\rho:\tilde{X}\to X$ and an elliptic fibration $\tilde{\pi}:\tilde{X}\to \B$ with $\tilde\pi=\rho\circ\pi$.  Over the generic point of $\ell$, $\tilde{X}$ has fibers of Kodaira type $I_2$.  The Hodge numbers of $\tilde{X}$ are $(3,231)$.

The elementary but crucial point is that since we have blown up a locus that lives over $\Delta$, the elliptic fibrations $\pi$ and $\tilde\pi$ are isomorphic over $U$:
\[
\begin{array}{ccc}
\tilde\pi^{-1}(U) &\stackrel{\simeq}{\rightarrow}&\pi^{-1}(U)\\
\tilde\pi\downarrow\phantom{\tilde\pi}&&\phantom{\pi}\downarrow\pi\\
U&=&U
\end{array},
\]
so $(R^1\tilde\pi_*\bR_{\tilde{X}})|_U=(R^1\pi_*IC(X))|_U=\bV$.  The conclusion is that $H^*(\B,IC(\bV))$ can be computed using either $X$ or $\tilde{X}$, and we will proceed using $\tilde{X}$.

This observation already clarifies a longstanding observation about
F-theory: that F-theory on the singular $X$ is dual to a limit of M-theory on its resolution $\tilde{X}$.

We claim that the Decomposition Theorem for $\tilde\pi$ reads
\begin{equation}\label{eq:decompI2}
R\tilde\pi_*\bR_{\tilde X}=\bR_B\oplus IC(\bV)[-1] \oplus \bR_\ell[-2]\oplus \bR_B[-2].
\end{equation}
As before, we confirm (\ref{eq:decompI2}) by restricting to a generic curve $C\subset \B$ and applying (\ref{eq:mdccurves}). We write $\Delta=\Delta_1\cup\ell$, the union of $\ell$ and an irreducible curve $\Delta_1$ of degree 34 obtained from (\ref{eq:i2disc}) after removing the factor of $\ell^2$.  The locus of critical values $\Sigma\subset C$ is $C\cap \Delta$.  Over points $p\in C\cap \Delta_1$, the fiber of $\tilde\pi$ is irreducible, hence the stalk $(T_\Sigma)_p$ vanishes.  Over points $p\in C\cap \ell$, the fiber of $\tilde\pi$ is an $I_2$ configuration, hence the stalk $(T_\Sigma)_p$ is 1-dimensional.

Reasoning as before using the Goresky-MacPherson inequality, we see that  $R\tilde\pi_*\bR_{\tilde X}$ can only have $\B$ and $\ell$ as supports, with the local systems associated with the corresponding IC sheaves being defined on $U$ and $\ell-Z$.  The local systems on $U$ can be found by restricting $R\tilde\pi_*\bR_{\tilde X}$ to $U$, and we get the summands $\bR_B\oplus IC(\bV)[-1] \oplus \bR_B[-2]$ exactly as before.  Furthermore, when we restrict $R\tilde\pi_*\bR_{\tilde X}$ to $C$, the torsion sheaf $T_\Sigma$ on $\Sigma=C\cap \ell$ has rank~1 at each $p\in \Sigma$, since $\tilde\pi^{-1}(p)$ is an $I_2$ configuration.  We conclude that the local system (call it $\bL$) on $\ell-Z$ has rank 1.  Finally, $\tilde\pi^{-1}(\ell)$ has two irreducible components: the exceptional divisor of $\rho$, and the proper transform of $\pi^{-1}(\Delta)$.  Recalling that the stalk of $T_\Sigma$ at $s\in\Sigma$ is $T_s\simeq H_2(f^{-1}(s))/\langle |f^{-1}(s)|\rangle$, we see that either of the two divisors trivializes $\bL$: $\bL\simeq\bR_{\ell-Z}$.  Thus the final contribution to (\ref{eq:decompI2}) comes from $IC(\bR_{\ell-Z})=\bR_\ell$.  Keeping track of the shifts required by our conventions, we have demonstrated (\ref{eq:decompI2}).

Taking $\bH^3$ of both sides of (\ref{eq:decompI2}) we get
\[
H^3(\tilde{X},\bR)\simeq\bH^3(\bP^2,\bR_{\bP^2}\oplus IC(\bV)[-1]\oplus \bR_\ell[-2]\oplus\bR_{\bP^2}[-2])
\]
\[=H^3(\bP^2,\bR)\oplus \bH^2(\bP^2,IC(\bV))\oplus H^1(\bP^2,\bR)\oplus H^1(\ell,\bR).
\]
Since all of the odd cohomologies on the right hand side vanish, this isomorphism simplifies to 
$\bH^2(\bP^2,IC(\bV))=H^3(\tilde{X},\bR)=\bR^{464}$, and we get 464 real scalars from dimensionally reducing the B-fields.  These scalars form half of the $h^{2,1}(\tilde{X})+1=232$ quaternionic scalars expected from F-theory.  These also match the number of uncharged scalars noted in Section~\ref{sec:examples}.

Taking $\bH^2$ of both sides of (\ref{eq:decompI2}) we get
\[
H^2(\tilde{X},\bR)=\bH^2(\bP^2,\bR_{\bP^2}\oplus IC(\bV)[-1]\oplus \bR_\ell[-2]\oplus\bR_{\bP^2}[-2])
\]
\[=H^2(\bP^2,\bR)\oplus \bH^1(\bP^2,IC(\bV))\oplus H^0(\bP^2,\bR)\oplus H^0(\ell,\bR),
\]
which implies that $\bH^1(\bP^2,IC(\bV))=0$, and no gauge fields arise from dimensionally reducing the B-fields.   This result matches the Shioda-Tate-Wazir formula perfectly assuming (\ref{eq:mwic}), as for the two components of $\Delta$ we have $n_1=1$ and $n_2=2$.

\smallskip
Having gone through these two examples, we can now see without much
more difficulty that our techniques go through more generally.
Express $\Delta=\Delta_1\cup\ldots\cup\Delta_k$ as the union of its
components.  Suppose that the fibers over a generic point of each
$\Delta_i$ have Kodaira type $I_n,\ I_n^*,\ II^*,\ III^*$, or $IV^*$,
and suppose further that we have a crepant resolution
$\rho:\tilde{X}\to X$ whose nontrivial fibers are ADE configurations
of $\bP^1$'s over generic points of the $\Delta_i$.  Letting
$\tilde\pi=\pi\circ\rho:\tilde{X}\to \B$, the fibers of $\tilde\pi$
over a generic point of $\Delta_i$  form
an affine ADE configuration of
curves.  For each $i$, let $U_i\subset \Delta_i$ be a Zariski open
subset contained in the smooth locus of $\Delta_i$ over which the
elliptic fibers of $\tilde\pi$ are equisingular.  We then have local
systems $\bL_i$ on $U_i$ whose stalk $(\bL_i)_p$ at $p\in U_i$ is
$H_2(\tilde\pi^{-1}(p))/\langle \tilde\pi^{-1}(p)\rangle$.  In this
situation, the decomposition theorem
reads
\begin{prop}\label{prop:decompgen}
\begin{equation*}
    R\tilde\pi_*\bR_{\tilde X}\simeq\bR_B\bigoplus IC(\bV)[-1] \bigoplus \oplus_i IC(\bL_i)[-2]\bigoplus \bR_B[-2].
\end{equation*}
\end{prop}
Note that the local systems $\bL_i$ can have monodromy, corresponding to the well-known phenomenon of monodromies in the exceptional curves of families of ADE resolutions.

\smallskip
Taking $\bH^2$ of both sides in the statement of Proposition~\ref{prop:decompgen}, we get
\begin{equation}\label{eq:matchstw}
    H^2(\tilde{X},\bR)\simeq H^2(\B,\bR) \bigoplus \bH^1(\B,IC(\bV)) \bigoplus\oplus_i \bH^0(\Delta_i,IC(\bL_i)) + 1
\end{equation}
We claim that $\dim \bH^0(\Delta_i,IC(\bL_i))=n_i-1$, so that (\ref{eq:matchstw}) matches the Shioda-Tate Wazir formula perfectly in complete generality, assuming (\ref{eq:mwic}).

Toward this end, we let $j_i:U_i\hookrightarrow \Delta_i$ be the inclusion.   We first claim that 
\begin{equation}\label{eq:icji}
    \bH^0(\Delta_i,IC(\bL_i))\simeq H^0(\Delta_i,(j_i)_*\bL_i).
\end{equation}
To see this, we first note that $\cH^0(IC(\bL_i))\simeq (j_i)_*\bL_i$, $\cH^1(IC(\bL_i))$ is a skycraper sheaf $T_i$ supported on the complement $Z_i$ of $U_i$ in $\Delta_i$, and $\cH^p(IC(\bL_i))=0$ for $p\ne0,1$.  These statements follow immediately from the general construction of $IC(\bL_i)$ in \cite[Section~4.2]{dmitm}.

It then follows that have the exact triangle in the constructible derived category of $\Delta_i$
\begin{equation}\label{eq:tri}
    (j_i)_*\bL_i \rightarrow IC(\bL_i) \rightarrow T_i[-1] \stackrel{+}{\rightarrow}.
\end{equation}
Then (\ref{eq:icji}) follows immediately from (\ref{eq:tri}) upon applying the long exact sequence of hypercohomology and noticing that $\bH^i(T_i[-1])=0$ unless $i=1$.

Let $\tilde\pi^{-1}({\Delta_i})$ be the union of irreducible divisors $D_{i1}\cup\ldots\cup D_{in_i}$.  
For each $1\le j\le n_i$ and $p\in U_i$, we have an element
$(s_j)_p$ of the stalk $(\bL_i)_p$ of $\bL_i$ at $p$ given by
\begin{equation}
    (s_j)_p=[\langle(\tilde\pi|_{D_{ij}})^{-1}(p)\rangle],
\end{equation}
the image of $\langle(\tilde\pi|_{D_{ij}})^{-1}(p)\rangle\in H_2(\tilde\pi^{-1}(p))$ in $(\bL_i)_p$.
The $(s_j)_p$  give a section $s_j\in H^0(\Delta_i,(j_i)_*\bL_i)$ by varying $p$, and 
\begin{equation}\label{eq:reln}
 \sum_{j=1}^{n_i}s_j=0   
\end{equation}
by the definition of $\bL_i$.
Each $(s_j)_p$ is the sum of the cohomology classes of the components of the affine ADE configuration of curves $\tilde\pi^{-1}(p)$ which lie in the irreducible component $D_{ij}$, and as a consequence this sum is monodromy invariant.  From this description we see that the $s_j$ span all monodromy invariant combinations, i.e.\ the $\{s_j\}$ span $H^0(\Delta_i,(j_i)_*\bL_i)$.  We also see easily that (\ref{eq:reln}) is the only relation among the $s_j$.  Hence $\dim H^0(\Delta_i,(j_i)_*\bL_i)=n_i-1$, so by (\ref{eq:icji}) we also have $\dim \bH^0(\Delta_i,IC(\bL_i))=n_i-1$, as claimed.

Comparing (\ref{eq:matchstw}) with the Shioda-Tate-Wazir formula (\ref{eq:stw}), we have proven

\begin{prop}
Assuming the hypotheses on $X\to B$ stated immediately before Proposition~\ref{prop:decompgen} together with  $H^{2,0}(X)=H^{2,0}(B)=0$, we have $\operatorname{rank}\mathrm{MW}(X)=\dim\bH^1(\B,IC(\bV))$.
\end{prop}

\medskip
The analysis of the scalar fields in more general cases requires a more detailed analysis of the 7-branes  \cite{kt}.

\section{Mordell-Weil} \label{sec:MW}

In the previous section, we saw that our conjectures and calculations show that 
\begin{equation}
    \operatorname{rank}(MW(X))=\dim \bH^1(\B,IC(\bV)),
\end{equation}
so that $MW(X)$ generates the gauge fields coming from dimensional reduction of the B-fields.  In this section, we make this equality more precise by showing that there is a group homomorphism
\begin{equation}\label{eq:mwmap}
    MW(X)\rightarrow \bH^1(\B,IC(\bV)).
\end{equation}
Such a map was constructed for elliptic surfaces (not necessarily K3 surfaces) in \cite{cz} using the Leray filtration.  

\subsection{A map from the Mordell-Weil group to intersection cohomology}\label{subsec:themap}

We adapt the idea of \cite{cz} by using the perverse Leray filtration in place of the usual Leray filtration.  

Our starting point is Proposition~\ref{prop:decompgen}, the Decomposition Theorem for $R\tilde\pi_*\bR_{\tilde X}$, which we use to describe a canonical map $R\tilde\pi_*\bR_{\tilde X}\to IC(\bV)[-1]$. While we can certainly use the existence of an isomorphism as in the statement of  Proposition~\ref{prop:decompgen}, and then project to the $IC(\bV)[-1]$ summand of the right hand side, the isomorphism is not canonical.  So we have to proceed with a bit more care.  

We can identify the perverse direct image sheaves $^{p} R^i\tilde\pi_*\bR_{\tilde X}$ by noting which shifts are needed to make each of the summands on the right hand side of Proposition~\ref{prop:decompgen} perverse.  The perverse shifts are
\begin{equation}
    ^\pi\bR_B=\bR_B[2],\ ^\pi\bV=\left(IC(\bV)[-1]\right)[3],\ ^\pi\bL_i=\left(IC(\bL_i)[-2]\right)[3].
\end{equation}
We conclude that
\begin{equation}\label{eq:pervcoh}
  ^pR^2\tilde\pi_*\bR_{\tilde X}\simeq \bR_B[2],\ 
  {}^pR^3\tilde\pi_*\bR_{\tilde X}\simeq IC(\bV)[2]\oplus IC(\bL_i)[1],\ 
  {}^pR^4\tilde\pi_*\bR_{\tilde X}\simeq \bR_B[2].
\end{equation}
Since the perverse sheaves $IC(\bV)[2]$ and $IC(\bL_i)[1]$ are simple, the projection map
\begin{equation}\label{eq:canproj}
    ^pR^3\tilde\pi_*\bR_{\tilde X}\rightarrow IC(\bV)[2]
\end{equation}
is canonical.

Denoting truncations with respect to the perverse t-stucture by $^p\tau$ as usual, (\ref{eq:pervcoh}) implies that 
$^p\tau_{\le4}(R\tilde\pi_*\bR_{\tilde X})=R\tilde\pi_*\bR_{\tilde X}$, and the canonical map
\begin{equation}
^p\tau_{\le4}(R\tilde\pi_*\bR_{\tilde X})\to {}^p\cH^4(R\tilde\pi_*\bR_{\tilde X})[-4]   
={}^pR^4\tilde\pi_*\bR_{\tilde X}[-4] 
\end{equation}
is identified with
$R\tilde\pi_*\bR_{\tilde X}\to \bR_B[-2]$.  Applying $\bH^2$ to $R\tilde\pi_*\bR_{\tilde X}\to \bR_B[-2]$, we get a canonical map
\begin{equation}
  \tilde\pi_*:  H^2(\tilde X,\bR)\to H^0(\B,\bR),
\end{equation}
which is identified with $\tilde\pi_*:H_4(\tilde X,\bR)\to H_4(\B,\bR)$ via Poincar\'e duality.  We also have the triangle
\begin{equation}
    ^p\tau_{\le3}(R\tilde\pi_*\bR_{\tilde X})\rightarrow {}^p\tau_{\le4}(R\tilde\pi_*\bR_{\tilde X})\rightarrow {}^p\cH^4(R\tilde\pi_*\bR_{\tilde X})[-4]\stackrel{+}{\rightarrow},
\end{equation}
which is identified with
\begin{equation}
    ^p\tau_{\le3}(R\tilde\pi_*\bR_{\tilde X})\rightarrow R\tilde\pi_*\bR_{\tilde X}\rightarrow\bR_B[-2]
    \stackrel{+}{\rightarrow},
\end{equation}
whose long exact hypercohomology sequence includes
\begin{equation}\label{eq:leh}
    0\rightarrow\bH^2(\tilde{X},{}^p\tau_{\le3}(R\tilde\pi_*\bR_{\tilde X}))\rightarrow H^2({\tilde X},\bR)\rightarrow H^0(\B,\bR_B).
\end{equation}

Now suppose that $s:\B --\to X$ is a rational section, which induces a rational section $\tilde{s}:\B\to \tilde X$.  The closure of $\tilde{s}$ is a divisor $D_{\tilde{s}}$ on $\tilde{X}$.  In particular, we have the divisor $D_0$ associated to the section of $X$ which gives $X$ the structure of an elliptic fibration.

Since $\tilde\pi_*([D_{\tilde{s}}]-[D_0])=0$, it follows from (\ref{eq:leh}) that $[D_{\tilde{s}}]-[D_0]\in H^2(\tilde{X},\bR)$ lives in the subspace $\bH^2(\tilde{X},{}^p\tau_{\le3}(R\tilde\pi_*\bR_{\tilde X}))$.  

The composition of the canonical maps 
\begin{equation}
{}^p\tau_{\le3}(R\tilde\pi_*\bR_{\tilde X})\to {}^p\cH^3(R\tilde\pi_*\bR_{\tilde X})[-3]={}^pR^3\tilde\pi_*\bR_{\tilde X}[-3]
\end{equation}
and the shift by $-3$ of ${}^pR^3\tilde\pi_*\bR_{\tilde X}\to IC(\bV)[2]$ (\ref{eq:canproj}) induces the map
\begin{equation}
   \sigma: \bH^2(\tilde{X},{}^p\tau_{\le3}(R\tilde\pi_*\bR_{\tilde X}))\rightarrow \bH^2(\tilde{X},IC(\bV)[-1])=\bH^1(\tilde{X},IC(\bV))
\end{equation}
on hypercohomology.  

We can finally define the desired map (\ref{eq:mwmap}) by sending $s$ to $\sigma([D_{\tilde{s}}]-[D_0])$.

\subsection{A map from the Mordell-Weil group to de Rham cohomology}

Given an element of the Mordell-Weil group, we write down an explicit 1-form on $U$ valued in $\bV$ corresponding to the section $s$ instead of exhibiting a class in $\bH^1(\B,IC(\bV))$.  We expect that this 1-form represents the image of $\sigma([D_{\tilde{s}}]-[D_0])$ under the restriction map
\begin{equation}
  \bH^1(\B,IC(\bV))  \rightarrow \bH^1(U,(IC(\bV))|_U) \simeq H^1(U,\bV)
\end{equation}
but we have not checked this.

Away from $\Delta$ and possibly finitely many additional points of $\B$ where $s$ is not defined, we can locally realize $s(b)$ as a point of $E_b=\bC/(\bZ+\bZ\tau(b))$, where $\tau(b)$ is subject to the usual monodromy transformation.  We write
\[
s(b) = f_1(b)\tau(b) + f_2(b), \qquad f_i(b)\in \bR/\bZ.
\]
We associate to $s$ the locally defined 1-forms
\[
(\omega_1,\omega_2)=(df_1,df_2).
\]
Note that the integer ambiguity of the $f_i$ disappears after taking the differential.  

The locally defined forms $df_1$ and $df_2$ are  real and $C^\infty$, but together they have an intrinsic holomorphic description in terms of the differential of $s$, which takes $T_b(\B)$ to $(T_\pi)_{s(b)}X$, the vertical tangent space of $X$ at $s(b)$.  This gives a section $ds$ of $s^*(T_\pi)\otimes \Omega^1_B$, interpreted as a vector bundle on $U$ minus finitely many points.  
By Hartog's Theorem, $ds$ extends to all of $U$, and so the differential forms $df_i$ obtained by expressing $ds$ in real coordinates extend uniquely to (locally defined) $C^\infty$ forms at any of the finitely many points of $U$ where $s$ is not regular.

We next check that $(\omega_1,\omega_2)$ transforms as a section of $\bV$.  Letting $\tau'=(a\tau+b)/(c\tau +d)$, recall that the isomorphism $\bC/(\bZ+\bZ\tau')\to \bC/(\bZ+\bZ\tau)$ is induced by multiplication by $c\tau +d$ on $\bC$.  This gives
\[
f_1'\tau'+f_2'\mapsto (c\tau+d)(f_1'\tau'+f_2')=(a\tau+b)f_1'+(c\tau+d)f_2',
\]
which implies
\[
f_1=af_1'+cf_2'\qquad f_2=bf_1'+df_2',
\]
and consequently
\begin{equation}\label{eq:mwmonodromy}
\omega_1=a\omega_1'+c\omega_2'\qquad \omega_2=b\omega_1'+d\omega_2',
\end{equation}
which is immediately confirmed as the monodromy transformation of $\bV$.  The forms $\omega_1,\omega_2$ are closed, so $(\omega_1,\omega_2)$ represents an element of $H^1(U,\bV)$, as claimed.

\smallskip
It has been noted in the physics literature that torsion elements of $MW(X)$ do not contribute gauge fields.  We now have a simple explanation for this.  If $s$ is torsion of order $n$, then locally we can write
\begin{equation}
    s(b)=\frac{a_1}n\tau+\frac{a_2}n
\end{equation}
with the $a_i$ integers.  So $f_i=a_i/n$ is locally constant, and $\omega_i=df_i=0$.

\section{The brane/bulk decomposition of fields}\label{subsec:supports}

In this section, we present evidence for the holography-inspired
decomposition of the supergravity fields proposed in
Section~\ref{subsec:normalization} into normalizable ``bulk'' fields
and non-normalizable ``brane'' fields,
giving a concrete instance of the proposed exact sequence
(\ref{eq:bulkboundary}) in the case of an elliptically fibered K3
surface, 8 dimensional F-theory, and exploring the possibility of
extending this structure to 6D F-theory models.

\subsection{Cohomology with supports}

We start by reviewing cohomology with supports, both its global version and local versions, along with the local to global spectral sequence that relates them. 

Let $X$ be a topological space, $Z\subset X$ a closed subset.  Put $U=X-Z$, and let $i:Z\hookrightarrow X$ and $j:U\hookrightarrow X$ be the inclusions.   In our application, $X$ will be the F-theory base $\B=\bP^1$ and $Z$ will be the discriminant divisor $\Delta$.

If $F$ is a sheaf (of abelian groups) on $X$, we define $H^0_Z(X,F)$ to be the kernel of the restriction map
\[
H^0(X,F)\to H^0(U,F|_U).
\]
In other words, $H^0_Z(X,F)$ is the group of sections of $F$ with support on $Z$.  By definition, we have a short exact sequence
\begin{equation}\label{eq:h0z}
0\to H^0_Z(X,F)\to H^0(X,F)\to H^0(U,F|_U)
\end{equation}
for any sheaf $F$.   The derived functors of $F\mapsto H^0_Z(X,F)$ are denoted by $H^i_Z(X,F)$ and are referred to as the cohomologies of $F$ with support on $Z$.  The exact sequence (\ref{eq:h0z}) extends to a long exact sequence
\begin{equation}\label{eq:supportles}
0\to H^0_Z(X,F)\to H^0(X,F)\to H^0(U,F|_U)\to H^1_Z(X,F)\to H^1(X,F)\to H^1(U,F|_U)\to \ldots .
\end{equation}

\smallskip
We next recall the local version of cohomology with supports.  Returning to the general case, let $\cH^0_Z(F)\subset F$  be the subsheaf of $F$ whose local sections are those sections of $F$ which are supported on $Z$.   We view the $\cH^0_Z(F)$ as sheaves on $Z$.  When we want to view these as sheaves on $X$ we will write them as $i_*\cH^0_Z(X,F)$. The derived functors of $F\mapsto \cH^0_Z(F)$ are denoted $\cH^i_Z(F)$, also understood as sheaves on $Z$.  

Then cohomologies with support can be computed using the spectral sequence 
\begin{equation}\label{eq:localglobal}
E_2^{p,q}=H^p(Z,\cH^q_Z(F))=> H^{p+q}_Z(X,F).    
\end{equation}

\subsection{The brane/bulk decomposition in 8-dimensional F-theory}

Now we return to the case of 8D F-theory: $X=\bP^1$, $Z=\Delta=\{r_1,\ldots,r_{24}\}$, and $F=j_*\bV$.  Then from (\ref{eq:supportles}) we get
\begin{equation}\label{eq:holseq}
H^1_\Delta(\bP^1,j_*\bV)\rightarrow H^1(\bP^1,j_*\bV)\rightarrow H^1(U,\bV)\rightarrow H^2_\Delta(\bP^1,j_*\bV)\rightarrow H^2(\bP^1,j_*\bV),
\end{equation}
since $(j_*\bV)|_U$ is simply $\bV$.

We make the following claims:

\smallskip\noindent
{\bf Claims.}  
\begin{enumerate}
    \item $H^2(\bP^1,j_*\bV)=0$.
    \item $H^1_\Delta(\bP^1,j_*\bV)=0$.
        \item $H^2_\Delta(\bP^1,j_*\bV)=H^0(\Delta,\cH^2_\Delta(\bP^1,j_*\bV))$.
    \item $\cH^2_\Delta(j_*\bV)$ is a skycraper sheaf on $\Delta$, with 1-dimensional stalks over the points of $\Delta$.
\end{enumerate}

These claims imply what we want.  The exact sequence (\ref{eq:holseq}) becomes
\begin{equation}\label{eq:holses}
0\rightarrow H^1(\bP^1,j_*\bV)\rightarrow H^1(U,\bV)\rightarrow H^0(\Delta,\cH^2_\Delta(\bP^1,j_*\bV))\rightarrow 0.
\end{equation}
Comparing to (\ref{eq:bulkboundary}), the first term in (\ref{eq:holses}) is the space of cohomology classes of normalizable 1-forms, the second term are all of the 1-forms, and the third term are 0-forms on $\Delta$ after trivializing the sheaf $\cH^2_\Delta(\bP^1,j_*\bV)$, and we have achieved our goal.

\smallskip
It remains to verify these claims.  The first claim was already demonstrated at the end of the discussion of the elliptically fibered K3 case in Section~\ref{sec:main-conjecture}.

We begin by describing the sheaves $\cH^p_\Delta(j_*\bV)$.  In the general situation of a sheaf $F$ on $U=X-Z$, we have an exact sequence of sheaves on $X$ \cite[Section 58.77]{stacks}
\begin{equation}\label{eq:localsupports}
0\to i_*\cH^0_Z(F)\to F\to j_*(F|_U)\to i_*\cH^1_Z(F)\to 0,    
\end{equation}
and for each $p>1$ an isomorphism $R^pj_*(F|_U)\simeq i_*\cH^{p+1}_Z(F)$.

We interpret (\ref{eq:localsupports}) and the isomorphism which followed it in the situation of 8-dimensional F-theory, where $X=\bP^1$ and $F=\bV$. Recalling that $(j_*\bV)|_U=\bV$, we also see that $j_*(j_*\bV|_U)$ is simply $j_*\bV$.  Then (\ref{eq:localsupports}) becomes an exact sequence
\begin{equation}\label{eq:locsup8d}
0\to i_*\cH^0_\Delta(j_*\bV)\to j_*\bV\to j_*\bV\to i_*\cH^1_\Delta(j_*\bV)\to 0.    
\end{equation}

Since the map $j_*\bV\to j_*\bV$ in (\ref{eq:locsup8d}) is the identity, we see that $\cH^0_\Delta(j_*\bV)=\cH^1_\Delta(j_*\bV)=0$.  
The spectral sequence (\ref{eq:localglobal}) then tells us that $H^0_\Delta(\bP^1,j_*\bV)=H^1_\Delta(\bP^1,j_*\bV)=0$ as well, and in particular we have proven the second claim.  

The spectral sequence and the vanishing of $\cH^p_\Delta(j_*\bV)$ for $p=0,1$  also tells us that $H^2_\Delta(\bP^1,j_*\bV)= H^0(\Delta,\cH^2_\Delta(j_*\bV))$, which proves the third claim.  

So to compute $H^2_\Delta(\bP^1,j_*\bV)$, we only need to compute the sheaf $\cH^2_\Delta(j_*\bV)$ and show that it is a skyscraper sheaf which is 1-dimensional at each point $r \in\Delta$, proving the final claim.

The isomorphisms $R^pj_*(F|_U)\simeq i_*\cH^{p+1}_Z(F)$ become $R^pj_*\bV\simeq i_*\cH^{p+1}_\Delta(j_*\bV)$ for $p>1$.  In particular $i_*\cH^2_\Delta(j_*\bV)\simeq R^1j_*\bV$.  So we just have to compute the stalks $(R^1j_*\bV)_{r_i}$. 

In general, for any map $f:X\to Y$ and sheaf $F$ on $X$ we can compute the stalk of $R^pf_*F$ at $y\in Y$ as
\begin{equation}
\left(R^pf_*F\right)_y=\lim_{p\in V}\Gamma(f^{-1}(V),F),    
\end{equation}
where the limit is taken over open sets $V\subset Y$ containing $p$.  Proceeding with the computation of $(R^1j_*\bV)_{r_i}$, we can take our open sets $V$ to be discs containing $r_i$, sufficiently small that they do not contain any other point of $\Delta$.  For all such $V$, we have $j^{-1}(V)=V-\{r_i\}$, and the local systems $\bV_{V-\{r_i\}}$ are topologically the same: the are all completely described by the monodromy transformation $T$.  The maps in the direct limit are then isomorphisms, and any $V$ can be used to compute the stalk.  We get
\[
\left(R^1j_*\bV\right)_{r_i}=H^1(\Delta^*,\bV),
\]
where we have identified $V-\{r_i\}$ with a punctured disc $\Delta^*$, and $\bV$ is the rank 2 local system with monodromy $T$.  We can describe $\bV$ in terms of a basis $\{e_1,e_2\}$ and 
\[
T(e_1)=e_1+e_2,\qquad T(e_2)=e_2.
\]
Noting that $e_2$ spans a trivial rank 1 local system $\bC\cdot e_2$ isomorphic to $\bC$ contained in $L$, and $T$ is trivial on $e_1$ after modding out $L$ by $\bC\cdot e_2$, we get a short exact sequence of local systems on $\Delta^*$
\[
0\to \bC\to \bV\to \bC\to 0.
\]
Computing cohomologies, we get
\begin{equation}
0\to H^0(\Delta^*,\bC)\rightarrow H^0(\Delta^*,\bV)\rightarrow H^0(\Delta^*,\bC)\rightarrow H^1(\Delta^*,\bC)\rightarrow H^1(\Delta^*,\bV)\rightarrow H^1(\Delta^*,\bC) \rightarrow 0.
\end{equation}

Using $H^0(\Delta^*,\bC)=H^1(\Delta^*,\bC)=H^0(\Delta^*,\bV)=\bC$, we deduce that $H^1(\Delta^*,\bV)\simeq\bC$, so that $R^1j_*\bV$ has the indicated structure as a skycraper sheaf, completing the verification of the last claim.

\subsection{The brane/bulk decomposition in 6-dimensional F-theory}

We cannot use the methods of the preceding section for analyzing the field decomposition in 6D F-theory, since $IC(\bV)$ can differ from $j_*\bV$ at the singularities of $\Delta$, as noted in the $n=2$ example in Section~\ref{subsec:mathexamples}.  Instead, we proceed using the Decomposition Theorem on both $X$ and $\pi^{-1}(U)$, and using the long exact sequence of Borel-Moore homology to compare them.  

The Decomposition Theorem for $\tilde{X}$ in the general 6D case is given by Proposition~\ref{prop:decompgen}.  Since $\tilde\pi|_U:\pi^{-1}(U)\to U$ is smooth, Deligne's Theorem applies just as it did in the 8D case and we get, putting $Y=\tilde\pi^{-1}(U)$

\begin{equation}\label{eq:6ddecompu}
\left(R\tilde\pi|_U\right)_*\bR_{Y}=\bR_U\oplus \bV[-1]\oplus \bR_U[-2].
\end{equation}

Then the restriction map $\rho:H^{p+1}(\tilde{X},\bR)\to H^{p+1}(Y,\bR)$ can be described by using the Decomposition Theorem on $\tilde{X}$ and $Y$ as just described, and applying hypercohomology.  The result is a map
\begin{equation}\label{eq:restdecomp}
\begin{split}
       H^{p+1}({\cal B},\bR)\bigoplus \bH^{p}({\cal B},IC(\bV))\bigoplus \oplus_i 
       \bH^{p-1}(\Delta_i,\bL_i)\bigoplus H^{p-1}({\cal B},\bR)\rightarrow\\  
       H^{p+1}(U,\bR)\bigoplus H^{p}(U,\bV)\bigoplus H^{p-1}(U,\bR).
\end{split}
\end{equation}
We see that the map $\bH^{p}({\cal B},IC(\bV))\to H^{p}(U,\bV)$ is a summand of $H^{p+1}(\tilde{X},\bR)\to H^{p+1}(Y,\bR)$.

We now put $D=(\tilde\pi)^{-1}(\Delta)$ and consider the long exact sequence of Borel-Moore homology
\begin{equation}\label{eq:BM}
 \cdots\rightarrow H^{BM}_{7-p}(Y)\rightarrow   H^{BM}_{6-p}(D)\rightarrow H^{BM}_{6-p}(\tilde{X})\rightarrow H^{BM}_{6-p}(Y)\rightarrow H^{BM}_{5-p}(D)\rightarrow \cdots
\end{equation}
Since $X$ and $Y$ are smooth, we have $H^i(\tilde{X})\simeq H^{BM}_{6-i}(\tilde{X})$ and $H^i(Y)\simeq H^{BM}_{6-i}(Y)$ by Poincar\'e duality.  Since $D$ is compact, we have $H^{BM}_{6-i}(D)\simeq H_{6-i}(D)$.  So (\ref{eq:BM}) becomes
\begin{equation}\label{eq:BM2}
   \cdots\rightarrow H^{p-1}(Y)\rightarrow   H_{6-p}(D)\rightarrow H^{p}(\tilde{X})\rightarrow H^{p}(Y)\rightarrow H_{5-p}(D)\rightarrow \cdots  
\end{equation}

Since $H_5(D)=0$ for dimension reasons, we see that
$H^1(\tilde{X},\bR)\to H^1(Y,\bR)$ is injective.  Therefore its
summand $\bH^{0}({\cal B},IC(\bV))\to H^{0}(U,\bV)$ is injective as
well, consistent with the physical expectations of Conjecture~\ref{c:physics}.

We need an extra step to analyze the injectivity of $H^2(\tilde{X},\bR)\to H^2(Y,\bR)$ using (\ref{eq:BM2}), since $H_4(D,\bR)$ is nonzero.

If there were a nonzero element $\omega$ of the kernel of $H^1(\tilde{X},IC(\bV))\to H^1(U,\bR)$, it would correspond to a nonzero element $\tilde\omega$ of the kernel of $H^2(\tilde{X},\bR)\to H^2(Y,\bR)$ by (\ref{eq:restdecomp}).  However, every element of this kernel comes from $H_4(D)$ by (\ref{eq:BM2}).  We try to derive a contradiction.

By the first equality in (\ref{eq:pervcoh}), we see that $^pR^2\tilde\pi_*\bR_{\tilde X}=\bR_B[2]$.  In other words, $^p\tau_{\le2}$ defines a canonical map $(\bR_B[2])[-2]\to R\tilde\pi_*\bR_{\tilde X}$ which fixes the summand $\bR_B$ in Proposition~\ref{prop:decompgen}.

Taking $\bH^2$, we get a canonically embedded summand
\begin{equation}\label{eq:summand}
    H^2({\cal B},\bR) \hookrightarrow H^2(\tilde{X},\bR),
\end{equation}
of the source of (\ref{eq:restdecomp}) with $p=1$.  This map
is identified with $\tilde\pi^*$.

If $D$ were irreducible, then $[D]$ generates the kernel of $H^2(\tilde{X},\bR)\to H^2(Y,\bR)$, so $\tilde\omega=[D]$ up to a scalar multiple.  Since $D=\tilde\pi^*([\Delta])$, we see that $[D]$ lies in the subspace $H^2({\cal B},\bR)$ of $H^2(\tilde{X},\bR)$ associated with $^p\tau_{\le2}$ identified in (\ref{eq:summand}).  In particular, $\tilde\omega$ projects to zero in the summand $H^1({\cal B},IC(\bV))$.  But the projection of $\tilde\omega$ is just $\omega$ by the definition of $\tilde\omega$, and we have proven that $H^1(\tilde{X},IC(\bV))\to H^1(U,\bR)$ is injective in this case.

We defer the proof in the general case to future work.  Each irreducible component $D_i$ of $D$ projects to some component $\Delta_i$ of $\Delta$ with fibers of real dimension 2.  So it seems reasonable to expect that
$[D_i]\in H^2(\tilde{X},\bR)$ lies in the summand 
\begin{equation}
 H^2({\cal B},\bR) \bigoplus \oplus_i  \bH^{0}(\Delta_i,\bL_i)
\end{equation}
of the source of (\ref{eq:restdecomp}) with $p=1$.  This would ensure that $[D_i]$ projects trivially to $\bH^1({\cal B},IC(\bV))$, and then the proof of injectivity goes through unchanged, since the $[D_i]$ generate the kernel of $H^2(\tilde{X},\bR)\to H^2(Y,\bR)$.

We also expect a similar argument to go through in the case $p=2$.  The kernel of $H^3(\tilde{X},\bR)\to H^3(Y,\bR)$ is generated by $H_3(D,\bR)$.  We expect that the image of $H_3(D,\bR)$ in $H^3(\tilde{X},\bR)$ lies in the summand 
\begin{equation}
 H^3({\cal B},\bR) \bigoplus \oplus_i  \bH^{1}(\Delta_i,\bL_i)
\end{equation}
of $H^3(\tilde{X},\bR)$,
which would prove that the projection of any such image to $\bH^2({\cal B},IC(\bV))$ would be zero.  The proof of injectivity would then go through as before.  We defer the details to future work.

\smallskip
In summary, we have given a mathematical proof of the injectivity in (\ref{eq:bulkboundary}) in certain cases after identifying $H^p_{norm}({\cal B}\backslash\Delta)$ with $\bH^p({\cal B},IC(\bV))$ and 
$H^p_{all}({\cal B}\backslash\Delta)$ with $H^p(U,\bV)$,
and have found a plausibility argument in the general case.

\section{Conclusions}

We have conjectured a mathematical description of the gauge fields and
scalars arising from the dimensional reduction of normalizable
B-fields in F-theory, and shown that the results are always consistent
with the expectations of physics in many large classes of
important cases.  We have
also proposed a bulk-boundary correspondence and provided supporting
mathematical evidence.

Another way to interpret our results is to ignore all discussion of
normalizability in  Conjectures~\ref{c:physics}
and \ref{c:main},
and simply assert
that our results demonstrate that
 $\bH^1(B,IC(\bV))$ describes the physically relevant gauge fields
arising from dimensional reduction of the B-fields, and that
$\bH^2(B,IC(\bV))$ describes the physically relevant scalars arising
from dimensional reduction of the B-fields.  From this viewpoint,
our conjectures provide
a possible physical explanation of
these assertions.

The comparison of our computation of scalars with physics is
necessarily limited since we have not discussed the scalars living on
7-branes in this paper.  Indeed, our analysis is complete only in the
special situation where there are no nonabelian gauge fields living on
the 7-branes, and correspondingly no additional
charged
scalars.  A more
general analysis including 
degrees of freedom arising from scalars
on the 7-branes will appear elsewhere \cite{kt}.

\end{document}